\documentclass[12pt,reqno]{amsart}
\usepackage{amsmath,amssymb,amscd,verbatim,color}
\usepackage{enumerate}
\makeatletter
\@namedef{subjclassname@2020}{%
	\textup{2020} Mathematics Subject Classification}
\makeatother
\usepackage[colorlinks=true]{hyperref}
\usepackage{amsfonts}
\usepackage[top=35mm, bottom=35mm, left=30mm, right=30mm]{geometry}

\usepackage{xcolor}
\theoremstyle{plain}
\newtheorem{thm}{Theorem}[section]
\newtheorem{cor}[thm]{Corollary}

\newtheorem{lem}[thm]{Lemma}
\newtheorem{prop}[thm]{Proposition}

\theoremstyle{definition}
\newtheorem{defn}[thm]{Definition}
\newtheorem{exam}[thm]{Example}
\newtheorem{rem}[thm]{Remark}
\newtheorem{claim}{Claim}

\numberwithin{equation}{section}

\newcommand{\mb}{\mathbb}

\newcommand{\ip}[2]{\left\langle #1,#2\right\rangle}
\newcommand{\norm}[1]{\left\lVert #1\right\rVert}
\newcommand{\abs}[1]{\left\lvert #1\right\rvert}
\newcommand{\Ran}{\operatorname{Ran}}

 \def \b{\beta}

\def\R{\mathbb{R}}

\allowdisplaybreaks

\begin{document}
	
	\title[Directional Ergodic Theorems, Ergodicity, and Weak mixing for $\mb Z^2$-actions]{Directional Ergodic Theorems, Ergodicity, and Weak mixing for $\mb Z^2$-actions}

\author[C. Liu,  K. Liu, T. Yu]{Chunlin Liu,  Kairan Liu and Tao Yu}

\address[C. Liu]{School of Mathematical Sciences, Dalian University of Technology, Dalian, 116024, P.R. China, and
	Institute of Mathematics, Polish Academy of Sciences, ul. \'sniadeckich 8, 00-656 Warszawa, Poland.}
\email{chunlinliu@mail.ustc.edu.cn}

\address[K. Liu]{College of Mathematics and Statistics, Key Laboratory of Nonlinear Analysis and its Applications (Chongqing University), Ministry of Education, Chongqing University, Chongqing, 401331, P.R. China}
\email{lkr111@cqu.edu.cn}

\address[T. Yu]{Department of Mathematics, Shantou University, Shantou, 515821, P.R. China}
\email{ytnuo@mail.ustc.edu.cn}

\subjclass[2020]{Primary  37A30; Secondary 37A25.}

\keywords{Ergodic Theorems; Ergodicity; Weak mixing}
\begin{abstract}
	In this paper, we study Birkhoff and von Neumann type ergodic theorems for $\mb Z^2$-measure preserving systems along irrational directions.
This leads us to a further investigation of directional ergodicity and weak mixing.
Additionally, we show that our framework (dynamical definitions) is equivalent to the spectral definitions of directional ergodicity and weak mixing introduced by Robinson, Rosenblatt and \c{S}ahin  \cite{RRS2023} for $\mb Z^2$-actions.

We note that in their recent paper \cite{RS2026}, Robinson and  \c{S}ahin established the equivalence between dynamical and spectral directional weak mixing for $\mb Z^2$-actions using an approach fundamentally different from the one employed in the present paper.
\end{abstract}

\maketitle

\section{Introduction}

  The mean ergodic theorem of von Neumann \cite{vonNe32} and the pointwise ergodic theorem of
	Birkhoff \cite{Bir31} are cornerstones of ergodic theory. In this paper we study directional
	analogues of these theorems for $\mb Z^2$-measure preserving systems, where the averages are
	taken over lattice points lying in a thin strip around a fixed line.

	The study of directional phenomena goes back to Milnor's introduction of directional entropy in the context of cellular automata \cite{Mil86,Mil88}. Continuity and further properties of directional entropy were subsequently studied by Park \cite{Par94,Pa}. Later, Kami\'nski-Park extended directional entropy to arbitrary $\mathbb Z^2$-actions on Lebesgue spaces \cite{KamPark99}. As
	research deepened, directional dynamics for $\mathbb Z^d$-actions has been investigated from several viewpoints (see e.g. \cite{BoyleLind97,Bro16,Dan17,JohnsonSahin15,RRS2023,RS2026}). 
Robinson, Rosenblatt and \c{S}ahin \cite{RRS2023} investigated directional ergodicity and weak mixing from a spectral-theoretic perspective, and their work provided valuable inspiration for the present paper.
	
	However, no directional version of ergodic average theorems has been established so far, despite its considerable importance. To be more precise, let $(X,\mathcal{B}_X,\mu,T)$ be a $\mathbb{Z}^2$-MPS. For every map $f\colon X\to\mathbb{R}$ and non-empty finite set $\Lambda\subseteq\mathbb{Z}^2$ we denote
	$$S(f,\Lambda)(x)=\sum_{(m,n)\in \Lambda}f(T^{(m,n)}x)\ \ \text{and}\ \ \mathbb{E}(f,\Lambda)(x)\colon=\frac{1}{\#\Lambda}S(f,\Lambda)(x)$$
	for every $x\in X$.  For every direction vector \footnote{When $\beta\in\mathbb Q$ the corresponding averages reduce to a more classical one-dimensional
		situation, and the results become essentially trivial, hence throughout we assume
		$\beta\notin\mathbb Q$; we only need to consider the case $\beta>0$, otherwise, we consider $S^{(m,n)}=T^{(m,-n)}$.} $\vec{v}=(1,\b)\in\mathbb{R}^2$ and
	$b\in (0,\infty)$, we put
	$$\Lambda^{\vec{v}}(b)=\left\{(m,n)\in\mathbb{Z}^2:\beta m-b/2\leq n\leq \beta m+b/2\right\},$$
	and $\Lambda_M^{\vec{v}}(b)=\Lambda^{\vec{v}}(b)\cap\big(\{0,1,\cdots,M-1\}\times\mathbb{Z}\big)$ for every $M\in\mathbb{N}$. We would like to investigate whether or not the ergodic averaging property holds along $\Lambda^{\vec{v}}(b)$.
	
	 Note that for fixed $b$ the sets $\Lambda^{\vec v}_M(b)$ are typically non-F\o lner, and
	therefore the classical multiparameter ergodic theorems (see e.g., \cite{Lin01}) do not apply directly. Nevertheless, our first main result establishes directional pointwise and mean ergodic theorems in this non-F\o lner setting; see Theorem~\ref{thm-1} and Theorem~\ref{D-von Ne-thm} for stronger variants.
	\begin{thm}\label{main-thm1}
		Let $(X,\mathcal{B}_X,\mu,T)$ be a $\mathbb{Z}^2$-MPS and let $\vec v=(1,\beta)$ with $\beta \not\in \mathbb{Q}$.
Then for every $b>0$ and every $f\in L^1(\mu)$ there exists
		$f_b^*\in L^1(\mu)$ such that
		\[
		\mathbb{E}(f,\Lambda_M^{\vec v}(b))\to f_b^*
		\quad\text{both $\mu$-a.e.\ and in $L^1(\mu)$ as }M\to\infty,
		\]
		and moreover $\int_X f_b^*\,d\mu=\int_X f\,d\mu$.
	\end{thm}
	
	It is worth mentioning that for a given $f\in L^1(\mu)$ if we define $\Phi_f\colon(0,\infty)\to L^1(\mu)$ by $\Phi_f(b)=f^*_b$, then it is continuous (see Proposition \ref{continuous on b}).
When $b\to 0^+$, the points in $\Lambda^{\vec{v}}(b)$  lie in increasingly narrow strips around the
line spanned by $\vec v$.
This leads to the question:
does the limit $\lim_{b\to 0^+}f_b^*$ exist in $L^1(\mu)$?  We give a positive answer to this question,
see Theorem~\ref{thm:zero-width-limit} for stronger variants.
	
\begin{thm}\label{main-thm2}
There exists a $T$-invariant sub-$\sigma$-algebra
$\mathcal{F}_{\vec v}\subseteq \mathcal{B}_X$, unique modulo $\mu$-null sets, such that
for every $f\in L^1(\mu)$,
$$
 \lim_{b\to0^+}
 \norm{f_b^*-\mathbb E_\mu(f\mid \mathcal{F}_{\vec v})}_1=0.
$$
\end{thm}
	
	With these convergence theorems established, we may proceed to define ergodicity and weak mixing along directions (see Definition \ref{def-v-ergodic}). Notice that the definition of directional ergodicity using invariant sets is inadequate (see Proposition \ref{prop:reason}),  which distinguishes the
directional setting from the standard one.
		
		In fact, ergodicity and weak mixing along directions are well-established notions. Liu-Xu introduced and characterized directional discrete spectrum in terms of directional sequence entropy and bounded complexity \cite{LX,LX2}. From the spectral point of view, Robinson-Rosenblatt-\c{S}ahin introduced notions of directional ergodicity and weak mixing via unit suspensions and studied their properties \cite{RRS2023}. More recently, Liu introduced a notion of
	directional weak mixing and proved a directional version of the Koopman--von Neumann
	spectral mixing theorem \cite{Liu}. In the present paper, we prove the equivalence of these definitions (Theorem \ref{thm:ve=>ves}, Theorem~\ref{thm:weak mixing via spectrum}).

After completing this paper, we shared our manuscript with Professors  Robinson and  \c{S}ahin. They kindly shared with us their recent work \cite{RS2026}, which provides valuable insights into directional ergodicity and weak mixing from a spectral-theoretic perspective. In particular, they established the equivalence between dynamical and spectral directional weak mixing for $\mathbb{Z}^2$-actions.
Their proof \cite{RS2026} is based on an approach fundamentally different from the one developed in the present paper.

	\medskip
	\noindent\textbf{The structure of the paper is as follows.}  In Section~\ref{Pre} we recall basic concepts and
	preliminary results. Section~\ref{sec:Birkhoff's ergodic theorem} proves directional versions
	of Birkhoff's and von Neumann's ergodic theorems.
In Section~\ref{sec: Der-Dwm} we introduce
	directional ergodicity and directional weak mixing and relate them to the corresponding global
	properties of the system. Section~\ref{sec:ergodicity and mixing via spectrum} establishes
	the equivalence between the spectral and ergodic average formulations.
	Section~\ref{Sec:examples} provides illustrative examples. In Section~\ref{sec:main-thm2} we give a proof of Theorem~\ref{main-thm2}.


\section{Preliminaries}\label{Pre}

In this paper, $\R$, $\mb Q$, $\mb Z$, and $\mb N$ stand for the sets of real numbers, rational numbers, integers, and natural numbers, respectively. For a Lebesgue space $(X,\mathcal{B}_X,\mu)$ and $1\le p\le\infty$, we denote by $L^p(X,\mathcal{B}_X,\mu)$ the usual $L^p$-space. Whenever the underlying space and $\sigma$-algebra are clear, we simply write $L^p(\mu)$ instead.

\subsection{Measure preserving systems}

Let $G$ be a second countable, locally compact abelian group with identity
$e_G$, and let $(X,\mathcal{B}_X,\mu)$ be a Lebesgue
space. A \textbf{measurable $G$-action} on $X$ is a jointly measurable map
$T\colon G\times X\to X$, written $T(g,x)=T^g x$, such that
\[
T^{e_G}x=x,\qquad T^g(T^h x)=T^{gh}x
\quad\text{for all }g,h\in G\text{ and }x\in X.
\]
In particular, each $T^g$ is a bimeasurable bijection with inverse
$T^{g^{-1}}$. We call the action \textbf{measure preserving} if
\[
\mu\bigl(T^gA\bigr)=\mu(A)
\quad\text{for every }g\in G\text{ and }A\in\mathcal{B}_X.
\]
In this case, $(X,\mathcal{B}_X,\mu,T)$ is called a
\textbf{$G$-measure-preserving system} (briefly, a \textbf{$G$-MPS}).
In this paper, the group $G$ will usually be $\mathbb{Z}^d$ or $\mathbb{R}^d$ for $d=1,2$.

Let $(X,\mathcal{B}_X,\mu,T)$ be a $G$-MPS.
We consider the product $G$-MPS $(X\times X,\mathcal{B}_X\times \mathcal{B}_X,\mu\times \mu,T\times T)$, and the product action $T\times T$ is given by
  $$(T\times T)^{g}(x_1,x_2)=(T^{g}x_1,T^{g}x_2)$$
for every $g\in G$ and $x_1,x_2\in X$.

Let $(X,\mathcal{B}_X,\mu,T)$ be a $G$-MPS. We say it is \textbf{ergodic}
if every $G$-invariant set $A\in\mathcal{B}_X$ (i.e. $T^gA\overset{\mu}{=}A$ for any $g\in G$) has $\mu(A)=0$ or $1$, and we say it is \textbf{weakly mixing} if the product system $(X\times X,\mathcal{B}_X\times \mathcal{B}_X,\mu\times \mu,T\times T)$ is ergodic.

\subsection{Topological dynamical systems}
Let $G$ be a second countable, locally compact abelian group with identity
$e_G$,
 and let $X$ be a compact metric space.
A \textbf{topological $G$-action}
on $X$ is a jointly continuous map
$T\colon G\times X\to X$, written $T(g,x)=T^g x$, such that
\[
T^{e_G}x=x,\qquad T^g(T^h x)=T^{gh}x
\quad\text{for all }g,h\in G\text{ and }x\in X.
\]
In particular, each $T^g$ is a homeomorphism with inverse
$T^{g^{-1}}$.
In this case, $(X,T)$ is called a
\textbf{$G$-topological dynamical system} (briefly, a \textbf{$G$-TDS}).

Let $(X,T)$ be a $G$-TDS. A probability measure $\mu$ on $X$ is called \textbf{invariant}
if for every $g\in G$ and every $A\in\mathcal{B}_X$ one has $\mu(T^{g}A)=\mu(A)$, where $\mathcal{B}_X$ is the collection of all Borel sets in $X$. By the Krylov-Bogolioubov theorem there exists an invariant measure $\mu$ on $X$,
thus $(X,\mathcal{B}_X,\mu,T)$ is a $G$-MPS.

\begin{lem}\label{0625}
Let $(X,\mathcal{B}_X,\mu,T)$ be a $\mathbb{Z}^2$-MPS,
then there is a $\mathbb{Z}^2$-TDS $(Y,S)$ and a measure-theoretic isomorphism
 $\pi:(X,\mathcal{B}_X,\mu,T)\rightarrow(Y,\mathcal{B}_Y,\nu,S)$, where $\nu$ is an $S$-invariant Borel probability measure and $\mathcal{B}_Y$ is the Borel $\sigma$-algebra of $Y$.
\end{lem}

\subsection{Eigenfunctions and Kronecker algebra}

Let $G$ be a second countable, locally compact abelian group. The \textbf{dual group} $\widehat{G}$ of a group $G$ is the set of characters of the group, namely the set of all continuous homomorphisms $\gamma\colon G\to S^1=\{z\in\mathbb{C}\colon\vert z\vert=1\}$ with pointwise multiplication.

Let $(X,\mathcal{B}_X,\mu,T)$ be a $G$-MPS, $T$ induces operators
$$U_{T^g}\colon L^0(\mu)\to  L^0(\mu),\ \ U_{T^g}f=f\circ T^g$$
for every $g\in G$. In particular, $U_{T^g}$ restricts to a unitary operator on
$L^2(\mu)$. We say $f\in L^2(\mu)$ is an \textbf{eigenfunction corresponding to eigenvalue} $\gamma\in \widehat{G}$ if $U_{T^g}f(x)=\gamma(g)f(x)$ for all $g\in G$ and $\mu$-a.e. $x\in X$.

Let $(X,\mathcal{B}_X,\mu,T)$ be a $G$-MPS. The \textbf{Kronecker $\sigma$-algebra} $\mathcal{K}_{\mu}$ is given by
$$\mathcal{K}_{\mu}=\big\{B\in\mathcal{B}_X\colon \overline{\{U_{T^g}1_B\colon g\in G\}}\ \text{is compact in}\ L^{2}(\mu)\big\}.$$
Let $\mathcal{H}_d$ denote the closure of the linear span of eigenfunctions in $L^2(\mu)$. One has
\begin{align}\label{kro-eigen}
	\mathcal{H}_d=L^2(X,\mathcal{K}_{\mu},\mu)\ \ \text{(see e.g. \cite[Theorem 1.14]{KL16})}.
\end{align}
We shall use the following standard characterization of weak mixing
(see \cite[p.~67]{Gla03}).
\begin{thm}\label{X-weakly mixing} Let $G$ be a second countable, locally compact abelian group, $(X,\mathcal{B}_X,\mu,T)$ be a $G$-MPS. Then the following are equivalent:
	\begin{itemize}
	\item[(1)] $(X,\mathcal{B}_X,\mu,T)$ is weakly mixing;
		\item [(2)] $\mathcal{K}_{\mu}=\{\emptyset,X\}$;
		\item [(3)] 1 is the only eigenvalue and the only eigenfunctions are constants.	 \end{itemize}
\end{thm}

\section{Directional Birkhoff's ergodic theorem}\label{sec:Birkhoff's ergodic theorem}

In this section, we present the proof of Theorem \ref{main-thm1}, mainly with the aid of the following corresponding skew product. In fact, we shall establish stronger statements
(Theorem~\ref{thm-1}, Theorem~\ref{D-von Ne-thm}).

\subsection{A corresponding skew-product representation with respect to $\vec{v}$}
Let $(X,\mathcal{B}_X,\mu,T)$ be a $\mathbb{Z}^2$-MPS, and
$\vec{v}=(1,\beta)\in\mathbb{R}^2$ be a direction vector. We write
$\mathbb{T}:=\mathbb{R}/\mathbb{Z}$. Let $m$ denote the Lebesgue measure on $\mathbb{T}$,
and $\mathcal{C}$ be the Borel $\sigma$-algebra of $\mathbb{T}$.
To study directional averages, we construct a corresponding skew product as follows. Let
   \begin{align}\label{skew-product-action}
   	\widetilde{X}\colon=X\times \mathbb{T},\quad\widetilde{\mu}\colon=\mu\times m, \ \text{and}\ \ \widetilde{T}(x,u)\colon=(T^{(1,\lfloor\beta+u\rfloor)}x,R_{\beta}(u))
   \end{align}
   for every $(x,u)\in\widetilde{X}$, where  $R_\beta(u):=u+\beta \pmod 1$ and $\lfloor a\rfloor$ is the largest integer no greater than $a$. Let $\widetilde{\mathcal{B}}$ be the $\sigma$-algebra generated by $\mathcal{B}_X$ and $\mathcal{C}$. Clearly $\widetilde{T}$ is a $\widetilde{\mathcal{B}}$-measurable isomorphism preserving $\widetilde{\mu}$.

 In the following, we will consistently use $(\widetilde{X}, \widetilde{\mathcal{B}}, \widetilde{\mu}, \widetilde{T})$ to denote the corresponding skew product system associated with the system $(X, \mathcal{B}_X, \mu, T)$ with respect to $\vec{v}$.
 We shall use the following uniform absolute-continuity lemma for the fibers of a nonnegative invariant function.

\begin{lem}\label{key-lem}
Let $(X,\mathcal{B}_X,\mu,T)$ be a $\mathbb{Z}^2$-MPS,
$\vec{v}=(1,\beta)\in\mathbb{R}^2$ be a direction vector with $\beta \not\in \mathbb{Q}$.
Let $g\in L^1(\widetilde{\mu})$ be nonnegative and $\widetilde{T}$-invariant.  Then for every $a>0$ there exist $t\in(0,1)$ and an $R_\beta$-invariant Borel set $I\subseteq\mathbb T$ with $m(I)=1$ such that
\[
\int_Y g(x,u)\,d\mu(x)<a
\]
for every $u\in I$ and every $Y\in\mathcal{B}_X$ with $\mu(Y)<t$.
\end{lem}
\begin{proof}
For $M>0$, define
\[
H_M(u):=\int_X g(x,u) 1_{\{g(x,u)>M\}}\,d\mu(x),
\qquad u\in\mathbb T.
\]
Since $g$ is $\widetilde T$-invariant and each $T^{(m,n)}$ preserves $\mu$, for $m$-a.e. $u\in\mathbb T$ we have
\begin{align*}
H_M(R_\beta u)
&=\int_X g\big(T^{(1,\lfloor\beta+u\rfloor)}x,R_\beta u\big)
   1_{\{g(T^{(1,\lfloor\beta+u\rfloor)}x,R_\beta u)>M\}}\,d\mu(x)\\
&=H_M(u).
\end{align*} As $\beta\notin\mathbb Q$, the rotation $R_\beta$ is ergodic; hence $H_M$ is constant $m$-a.e. Denote this constant by $c_M$. By Tonelli's theorem,
\[
c_M=\int_{\mathbb T}H_M\,dm
=\int_{\widetilde X}g 1_{\{g>M\}}\,d\widetilde\mu
\longrightarrow0
\qquad\text{as }M\to\infty.
\]
Choose $M>0$ such that $c_M<a/2$, and set
\[
J:=\{u\in\mathbb T:H_M(u)=c_M\},
\qquad
I:=\bigcap_{n\in\mathbb Z}R_\beta^nJ.
\]
Then $I$ is an $R_\beta$-invariant Borel set with $m(I)=1$. Finally, choose
\[
0<t<\min\left\{1,\frac{a}{2M}\right\}.
\]
For every $u\in I$ and every $Y\in\mathcal{B}_X$ with $\mu(Y)<t$, we have
\begin{align*}
\int_Y g(x,u)\,d\mu(x)
&\leq \int_Y g(x,u) 1_{\{g(x,u)\leq M\}}\,d\mu(x)+H_M(u)\\
&\leq M\mu(Y)+c_M<a.
\end{align*}
This proves the lemma.
\end{proof}
\subsection{Directional ergodic average theorem}

For $\vec v=(1,\beta)\in \R^2$, $c_1<c_2\in\mathbb{R}$ and $M\in\mathbb{N}$, we define the lattice strip of width $c_2-c_1$ along the direction $\vec v$,
$$\Lambda_M^{\vec v}(c_1,c_2):=\bigl\{(m,n)\in \mathbb{Z}^2:\ \beta m+c_1\le n\le \beta m+c_2,\ 0\le m\le M-1\bigr\}.$$
The following theorem establishes the almost everywhere convergence of the ergodic averages over such lattice strips for every $f\in L^1(\mu)$.
 \begin{thm}\label{thm-1}
	Let $(X,\mathcal{B}_X,\mu,T)$ be a $\mathbb{Z}^2$-MPS,
$\vec{v}=(1,\beta)\in\mathbb{R}^2$ be a direction vector with $\beta \not\in \mathbb{Q}$.
Then for every $c_1,c_2\in\mathbb{R}$ with $c_1<c_2$ and $f\in L^1(\mu)$ there is a $f^*_{c_1,c_2}\in L^1(\mu)$ satisfying
	\begin{align}\label{Brik-eq}
	\lim_{M\to\infty}\mathbb{E}(f,\Lambda_M^{\vec{v}}(c_1,c_2))(x)=f^*_{c_1,c_2}(x)
	\end{align}
	for $\mu$-a.e. $x\in X$. Moreover, $\int_X f^{*}_{c_1,c_2}d\mu=\int_X fd\mu$.
\end{thm}
\begin{proof} By decomposing $f$ into its positive and negative parts we can assume $f\geq 0$. By decomposing the strip $[c_1,c_2]$ into finitely many subintervals of length less than 1, and by the $T$-invariance of $\mu$, we can also assume
$0<c_2-c_1<1$ and $0\le c_2<1$.

If $\int_X f\,d\mu=0$, then $f=0$ $\mu$-a.e. and there is nothing to prove. Thus, in what follows,
we assume $\int_X f\,d\mu>0$.
	
We set $\omega=c_2-c_1$. Let $(\widetilde{X},\widetilde{\mathcal{B}},\widetilde{\mu},\widetilde{T})$ be the corresponding skew product system associated with $(X,\mathcal{B}_X,\mu,T)$ and $\vec{v}$ (see (\ref{skew-product-action})). We define $\widetilde{f}\in L^{1}(\widetilde{\mu})$ by $$\widetilde{f}(x,u)=f(x)1_{[0, \ \omega]}(u)\ \text{for every}\ (x,u)\in \widetilde{X}.$$ Then by Birkhoff pointwise ergodic theorem, there exists a $\widetilde{T}$-invariant $\widetilde{f}^{*}\in L^1(\widetilde{\mu})$ and a $\widetilde{T}$-invariant subset $\widetilde{X}_0\subseteq \widetilde{X}$ with $\widetilde{\mu}(\widetilde{X}_0)=1$, such that for every $(x,u)\in\widetilde{X}_0$
	\begin{align}\label{Brikhoff-T-skew}
		\lim_{M\to \infty}\frac{1}{M}\sum_{i=0}^{M-1}\widetilde{f}(\widetilde{T}^i(x,u))=\widetilde{f}^*(x,u)<\infty,\  \text{and}\ \int_{\widetilde{X}} \widetilde{f}^*d\widetilde{\mu}=\int_{\widetilde{X}} \widetilde{f}d\widetilde{\mu}.
	\end{align}
	For every $u\in\mathbb{T}$, we set $F(u)\colon=\int_X \widetilde{f}^*(x,u)d\mu(x)$.
By Tonelli's theorem, $F$ is measurable, and
	$$\int_{\mathbb{T}} Fdm=\int_{\widetilde{X}} \widetilde{f}^*d\widetilde{\mu}=\int_{\widetilde{X}}\widetilde{f}d\widetilde{\mu}=\omega\int_X fd\mu.$$
	Moreover, note that $\mu$ is $T$-invariant and $\widetilde{f}^*$ is $\widetilde{T}$-invariant, for $m$-a.e. $u\in\mathbb{T}$, $$F(R_{\beta}(u))=\int_X \widetilde{f}^*(T^{(1,\lfloor\beta+u\rfloor)}x,R_{\beta}(u))d\mu(x)=\int_X \widetilde{f}^*(\widetilde{T}(x,u))d\mu(x)=F(u).$$
Since $\beta\notin\mathbb{Q}$, $R_{\beta}$ is ergodic and there exists a $R_{\beta}$-invariant set $I_0\subseteq\mathbb{T}$,
such that for every $u\in I_0$ one has
\begin{equation}\label{eq:F-constant-thm1}
	F(u)=\int_{\mathbb{T}}Fdm=\omega\int_X fd\mu.
\end{equation}

	Now we shall prove that for every $u_0\in\mathbb{T}$,
there exists a measurable $X_{u_0}\subseteq X$ with $\mu(X_{u_0})=1$ such that for every $x\in X_{u_0}$, $$\bigg\{\frac{1}{M}\sum_{i=0}^{M-1}\widetilde{f}\big(\widetilde{T}^i(x,u_0)\big)\bigg\}_{M\in\mathbb{N}}$$ is a Cauchy sequence.

We first suppose $u_0\in(0,1)$. For each $k\in\mathbb N$, apply Lemma~\ref{key-lem} to
$\widetilde f^{\,*}$ with $a=1/k$. Thus there are $\tau_k>0$ and an
$R_\beta$-invariant Borel set $I_k\subseteq\mathbb T$, $m(I_k)=1$, such that
\[
	\mu(Y)< \tau_k,\quad u\in I_k
	\quad\Longrightarrow\quad
	\int_Y\widetilde f^{\,*}(x,u)\,d\mu(x)< \frac1k.
\]
Shrinking the positive numbers if necessary, choose $t_k$ so that
\begin{equation}\label{eq:tk-choice-thm1}
	0<t_k<\min\left\{
	\frac1{k},\tau_k,
	\sqrt{u_0},\sqrt{\omega},\sqrt{1-\omega}
	\right\}.
\end{equation}
Then, for every $Y\in\mathcal B_X$ with $\mu(Y)\le t_k$ and every $u\in I_k$,
\begin{align}\label{e-thm-1}
	\int_Y\widetilde f^{\,*}(x,u)\,d\mu(x)<\frac1k.
\end{align}

We set $I_0'=\bigcap_{k\in\mathbb{N}}I_k$. Then $R_{\beta}I_0'=I_0'$ and $m(I_0')=1$. For every $k\in\mathbb{N}$, we define
	$$ \widetilde{g}_k(x,u)=f(T^{(0,1)}x)1_{[1-t_k^2,\ 1)}(u)+f(x)1_{[\omega-t_k^2,\ \omega]}(u)$$
	for every $(x,u)\in\widetilde{X}$. Then for every $k\in\mathbb{N}$ there exists a $\widetilde{T}$-invariant set $\widetilde{X}_k\subseteq\widetilde{X}_0$ with  $\widetilde{\mu}(\widetilde{X}_k)=1$,
and $\widetilde g^*_k\in L^{1}(\widetilde{\mu})$
with $\int_{\widetilde{X}} \widetilde{g}_kd\widetilde{\mu}=\int_{\widetilde{X}} \widetilde{g}^*_kd\widetilde{\mu}$,
 such that for every $(x,u)\in \widetilde{X}_k$
 one has $\widetilde g^*_k(x,u)=\widetilde g^*_k\big(\widetilde{T}(x,u)\big)$ and
	\begin{align}\label{Xk}
		\lim_{M\to \infty}\frac{1}{M}\sum_{i=0}^{M-1}\widetilde{g}_k(\widetilde{T}^i(x,u))=\widetilde{g}^*_k(x,u).
	\end{align}

	To show $\{\frac{1}{M}\sum_{i=0}^{M-1}\widetilde{f}(\widetilde{T}^i(x,u_0))\}_{M\in\mathbb{N}}$ is a Cauchy sequence for $\mu$-a.e $x\in X$, we need the following claims.
	 \begin{claim}\label{gk-f} For every $k\in\mathbb{N}$ and $u\in(u_0-t_k^2,u_0)$, for any   $M\in\mathbb{N}$ we have
		\begin{align}\label{eq:fiber-comparison-thm1}
\bigg\vert\sum_{i=0}^{M-1}\widetilde{f}(\widetilde{T}^i(x,u_0))
-\sum_{i=0}^{M-1}\widetilde{f}(\widetilde{T}^i(x,u))\bigg\vert
\leq\sum_{i=0}^{M-1}\widetilde{g}_k(\widetilde{T}^i(x,u)).
\end{align}
	\end{claim}
	\begin{proof}[Proof of Claim \ref{gk-f}] First, by the definition of $\widetilde{f}$, it is not hard to see
		\begin{align}\label{e-guji1}
			&\bigg\vert\sum_{i=0}^{M-1}\widetilde{f}(\widetilde{T}^i(x,u_0))
-\sum_{i=0}^{M-1}\widetilde{f}(\widetilde{T}^i(x,u))\bigg\vert\nonumber	\\
			=&\bigg\vert \sum_{i\in A_M\cup B_M}f(T^{(i,\lfloor i\beta+u_0\rfloor)}x)-\sum_{i\in A_M\cup B_M'}f(T^{(i,\lfloor i\beta+u\rfloor)}x)\bigg\vert,
		\end{align}
	where
	$$A_M=\big\{0\leq i\leq M-1\colon R_{\beta}^i(u_0),\ R_{\beta}^i(u)\in[0,\omega]\big\},$$
	$$B_M=\big\{0\leq i\leq M-1\colon R_{\beta}^i(u_0)\in[0,\omega]\ \text{and}\ R_{\beta}^i(u)\notin[0,\omega]\big\},$$
	$$B_M'=\big\{0\leq i\leq M-1\colon R_{\beta}^i(u_0)\notin[0,\omega]\ \text{and}\ R_{\beta}^i(u)\in[0,\omega]\big\}.$$
	We note that, since $u\in(u_0-t_k^2,u_0)$ and $t_k<\min\{\sqrt{1-\omega},\sqrt{\omega}\}$,
$$0<(i\beta+u_0)-(i\beta+u)=u_0-u<t_k^2<\min\{1-\omega,\omega\}.$$
For every $i\in A_M$ one has
$\lfloor i\beta+u_0\rfloor=\lfloor i\beta+u\rfloor$,
for every $i\in B_M$ one has $R_{\beta}^i(u)\in[1-t_k^2,1)$ and $\lfloor i\beta+u\rfloor+1=\lfloor i\beta+u_0\rfloor$,
for every $i\in B_M'$ one has $R_{\beta}^i(u)\in [\omega-t_k^2,\omega]$ and $\lfloor i\beta+u\rfloor=\lfloor i\beta+u_0\rfloor$.
Note that $f\geq 0$. Thus
		\begin{align*}
			(\ref{e-guji1})\leq & \sum_{i\in B_M}f(T^{(i,\lfloor i\beta+u\rfloor+1)}x)
+\sum_{i\in B_M'}f(T^{(i,\lfloor i\beta+u\rfloor)}x)\nonumber\\
			\leq & \sum_{i=0}^{M-1}\bigg(f(T^{(i,\lfloor i\beta+u\rfloor+1)}x)1_{[1-t_k^2,1)}(R_{\beta}^i(u))+f(T^{(i,\lfloor i\beta+u\rfloor)}x)1_{[\omega-t_k^2,\omega]}(R_{\beta}^i(u))\bigg)\nonumber\\
			=&\sum_{i=0}^{M-1}\widetilde{g}_k(\widetilde{T}^i(x,u))\nonumber.
		\end{align*}	
	\end{proof}	
	\begin{claim}\label{thm1-claim1} For every $k\in\mathbb{N}$,
there exists $u_k\in (u_0-t_k^2,u_0)\cap I_0\cap I_0'$, and $H_k\subseteq X$ with $\mu(H_k)\geq1-t_k$,
such that $H_k\times\{u_k\}\subset\widetilde{X}_k$
and for every $x\in H_k$ one has $\widetilde{g}^*_k(x,u_k)<2t_k\int_X fd\mu$.
	\end{claim}
	\begin{proof}[Proof of Claim \ref{thm1-claim1}] For every $k\in\mathbb{N}$, we set  $\widetilde{F}_k\colon=\{(x,u)\in\widetilde{X}_k\colon \widetilde{g}_k^*(x,u)<2t_k\int_X fd\mu\}$. It is $\widetilde{T}$-invariant. Then
		$$2t_k^2\int_X fd\mu=\int_{\widetilde{X}} \widetilde{g}_kd\widetilde{\mu}=\int_{\widetilde{X}} \widetilde{g}^*_kd\widetilde{\mu}\geq \int_{\widetilde{X}\setminus\widetilde{F}_k}\widetilde g^*_k(x,u)d\widetilde{\mu}
\geq\widetilde{\mu}(\widetilde{X}\setminus\widetilde{F}_k)\cdot2t_k\int_X fd\mu$$
		Thus $\widetilde{\mu}(\widetilde{F}_k)\geq 1-t_k$. Then there exists $u_k'\in I_0\cap I_0'$ and $H_k'\subset X$ with $\mu(H_k')\geq 1-t_k$ such that $H_k'\times\{u_k'\}\subseteq \widetilde{F}_k$. Since $\beta\notin \mathbb{Q}$, there exists $n_k\in\mathbb{N}$ such that $u_k\colon=R_{\beta}^{n_k}(u'_k)\in(u_0-t^2_k,u_0)$. Set $H_k\colon=T^{(n_k,\lfloor n_k\beta+u_k'\rfloor)}H_k'$. Then $\mu(H_k)\geq 1-t_k$
and $H_k\times\{u_k\}=\widetilde{T}^{n_k}(H_k'\times\{u_k'\})\subseteq \widetilde{F}_k$.
	\end{proof}
 Let $H=\bigcap_{n\in\mathbb{N}}\bigcup_{k\geq n}H_k$. Clearly $\mu(H)=1$. We claim that
  \begin{claim}\label{thm1-claim2}
  	For every $x\in H$ the following sequence is a Cauchy sequence:
  	\begin{align}\label{cauchyseq}
  	\bigg\{\frac{1}{M}\sum_{i=0}^{M-1}\widetilde{f}\big(\widetilde{T}^i(x,u_0)\big)\bigg\}_{M\in\mathbb{N}}.
  	\end{align}
  \end{claim}
  \begin{proof}[Proof of Claim \ref{thm1-claim2}]
   Fix $x\in H$. Given $\varepsilon>0$, take $k\in \mathbb{N}$ large enough such that $x\in H_k$ and $2t_k\int_X fd\mu<\varepsilon/5$. Since $(x,u_k)\in \widetilde{X}_k\subseteq\widetilde{X}_0$, by (\ref{Brikhoff-T-skew}) and (\ref{Xk}), there exists $M_0\in\mathbb{N}$ such that
  \begin{itemize}
  	\item [(i)] for every $M_1,M_2\geq M_0$ one has
     $$\bigg\vert\frac{1}{M_1}\sum_{i=0}^{M_1-1}\widetilde{f}(\widetilde{T}^i(x,u_k))
     -\frac{1}{M_2}\sum_{i=0}^{M_2-1}\widetilde{f}(\widetilde{T}^i(x,u_k))\bigg\vert<\frac{\varepsilon}{5},$$
     \item[(ii)] and for every $M\geq M_0$ one has
     $$\bigg\vert\frac{1}{M}\sum_{i=0}^{M-1}\widetilde{g}_k(\widetilde{T}^i(x,u_k))
     -\widetilde{g}^*_k(x,u_k)\bigg\vert<\frac{\varepsilon}{5}.$$
  \end{itemize}

  Then for every $M\geq M_0$, by Claim \ref{gk-f}, one has
  \begin{align}\label{eguji-2}
  	& \frac{1}{M}\bigg\vert\sum_{i=0}^{M-1}\widetilde{f}(\widetilde{T}^i(x,u_0))
  -\sum_{i=0}^{M-1}\widetilde{f}(\widetilde{T}^i(x,u_k))\bigg\vert\leq \frac{1}{M}\sum_{i=0}^{M-1}\widetilde{g}_k(\widetilde{T}^i(x,u_k))\nonumber\\
  	\overset{\text{(ii)}}{\leq}& \widetilde{g}^*_k(x,u_k)+\frac{\varepsilon}{5}
  \overset{(x\in H_k)}{\leq} 2t_k\int_X fd\mu +\frac{\varepsilon}{5}<\frac{2\varepsilon}{5}.
  \end{align}
  Hence for every $M_1,M_2\geq M_0$ we have
  \begin{align*}
  	&\bigg\vert\frac{1}{M_1}\sum_{i=0}^{M_1-1}\widetilde{f}(\widetilde{T}^i(x,u_0))
  -\frac{1}{M_2}\sum_{i=0}^{M_2-1}\widetilde{f}(\widetilde{T}^i(x,u_0))\bigg\vert\\
  	\overset{\text{(i)}}{\leq} &\frac{1}{M_1}\bigg\vert\sum_{i=0}^{M_1-1}\widetilde{f}(\widetilde{T}^i(x,u_0))
  -\sum_{i=0}^{M_1-1}\widetilde{f}(\widetilde{T}^i(x,u_k))\bigg\vert\\
  	&\ \ \ +\frac{1}{M_2}\bigg\vert\sum_{i=0}^{M_2-1}\widetilde{f}(\widetilde{T}^i(x,u_0))
  -\sum_{i=0}^{M_2-1}\widetilde{f}(\widetilde{T}^i(x,u_k))\bigg\vert+\frac{\varepsilon}{5}\\
  	\overset{(\ref{eguji-2})}{<} &\frac{4\varepsilon}{5}+\frac{\varepsilon}{5}=\varepsilon.
  \end{align*}
  Therefore, $(\ref{cauchyseq})$ is a Cauchy sequence.
  \end{proof}
 Thus there exists $\widehat f^*(x,u_0)\in\mathbb{R}$ such that for every $x\in H$
  \begin{align}\label{thm-con}
  	\frac{1}{M}\sum_{i=0}^{M-1}\widetilde{f}\big(\widetilde{T}^i(x,u_0)\big)\to \widehat f^{*}(x,u_0)\ \text{as}\ M\to\infty.
  \end{align}
 We claim that
 \begin{claim}\label{thm1-claim3}
 	$\int_{X}\widehat f^*(x,u_0)d\mu(x)=\omega\cdot \int_Xfd\mu(x).$
 \end{claim}
 \begin{proof}[Proof of Claim \ref{thm1-claim3}]
For $x\in H_k\cap H$, divide
\eqref{eq:fiber-comparison-thm1} by $M$ and let $M\to\infty$. Using
\eqref{Brikhoff-T-skew}, \eqref{Xk}, and Claim \ref{thm1-claim1}, we obtain
\[
	|\widehat f^{\,*}(x,u_0)-\widetilde f^{\,*}(x,u_k)|
	\le \widetilde g_k^{\,*}(x,u_k)<2t_k\int_X fd\mu.
\]
Therefore, by \eqref{eq:F-constant-thm1} and \eqref{e-thm-1},
\begin{align*}
	\int_X \widehat f^{\,*}(x,u_0)\,d\mu(x) &\geq \int_{H_k\cap H}\widehat f^*(x,u_0)d\mu(x)
	\ge
	\int_{H_k\cap H}\widetilde f^{\,*}(x,u_k)\,d\mu(x)-2t_k\int_X fd\mu\\
	&\ge \int_X\widetilde f^{\,*}(x,u_k)\,d\mu(x)-\frac{1}{k}-2t_k\int_X fd\mu \\
&= \omega \int_X fd\mu-\frac{1}{k}-2t_k\int_X fd\mu.
\end{align*}
 Let $k\to\infty$, then
\begin{equation}\label{eint1}
	\int_X \widehat f^{\,*}(x,u_0)\,d\mu(x)\ge\omega \int_X fd\mu.
\end{equation}

On the other hand, by Fatou's lemma and $\beta\notin \mathbb{Q}$,
  \begin{align}\label{eint2}
  	\int_X \widehat f^*(x,u_0)d\mu&\overset{(\ref{thm-con})}{\leq} \liminf_{M\to\infty} \frac{1}{M}\sum_{i=0}^{M-1}\int_X \widetilde{f}(\widetilde{T}^i(x,u_0))d\mu(x)\nonumber\\
  	&=\liminf_{M\to\infty} \frac{1}{M}\sum_{i=0}^{M-1}\int_X \widetilde{f}\big(T^{(i,\lfloor i\b+u_0\rfloor)}x,R_{\b}^i(u_0)\big)d\mu(x)\nonumber\\
  	&=\liminf_{M\to\infty} \frac{1}{M}\sum_{i=0}^{M-1}(\int_Xfd\mu(x))\cdot 1_{[0,\omega]}(R_{\b}^i(u_0))
  	=\omega\cdot \int_Xfd\mu(x).
  \end{align}
 By (\ref{eint1}) and (\ref{eint2}), we prove Claim \ref{thm1-claim3}.
 \end{proof}

It remains to treat $u_0=0$. Put $u_1=R_\beta(0)\in(0,1)$
and $\widetilde T(x,0)=(T^{(1,\lfloor\beta\rfloor)}x,u_1)$.
The preceding argument gives a full-measure set
$H_{u_1}$ on which
  $$\frac{1}{M}\sum_{i=0}^{M-1}\widetilde{f}\big(\widetilde{T}^i(x,u_1)\big)\to \widehat f^{*}(x,u_1)\ \text{as}\ M\to\infty.$$

For $\mu$-a.e. $x$ with $T^{(1,\lfloor\beta\rfloor)}x\in H_{u_1}$,
\[
\frac1M\sum_{i=0}^{M-1}\widetilde f(\widetilde T^i(x,0))
=
\frac{\widetilde f(x,0)}M
+
\frac{M-1}{M}
\left(
\frac1{M-1}\sum_{i=0}^{M-2}\widetilde f(\widetilde T^i(T^{(1,\lfloor\beta\rfloor)}x,u_1))
\right).
\]

Hence $\frac1M\sum_{i=0}^{M-1}\widetilde f(\widetilde T^i(x,0)) \rightarrow \widehat{f}^{\,*}(T^{(1,\lfloor\beta\rfloor)}x,u_1)$ as $M\to\infty.$
Denote $\widehat f^{*}(x,0):=\widehat{f}^{\,*}(T^{(1,\lfloor\beta\rfloor)}x,u_1)$.
Since $T^{(1,\lfloor\beta\rfloor)}$ preserves $\mu$,
$$\int_{X}\widehat f^{\,*}(x,0)d\mu(x)
=\int_{X}\widehat f^{\,*}(T^{(1,\lfloor\beta\rfloor)}x,u_1)d\mu(x)=\int_{X}\widehat f^{\,*}(x,u_1)d\mu(x)=
\omega\cdot \int_Xfd\mu(x).$$  Consequently, \eqref{thm-con} and
Claim \ref{thm1-claim3} hold for every prescribed $u_0\in[0,1)$.

For every $m\ge0$,
\[
	\widetilde T^m(x,u_0)
	=
	\left(T^{(m,\lfloor m\beta+u_0\rfloor)}x,
	R_\beta^m(u_0)\right).
\]
Since the interval
\[
	[m\beta+u_0-\omega,\,m\beta+u_0]
\]
has length $\omega<1$, it contains an integer if and only if
$R_\beta^m(u_0)\in[0,\omega]$; in that case the unique integer is
$\lfloor m\beta+u_0\rfloor$. Therefore, for every $M$,
\begin{equation}\label{eq:strip-skew-identity-thm1}
	S(f,\Lambda_M^{\vec v}(u_0-\omega,u_0))(x)
	=
	\sum_{m=0}^{M-1}\widetilde f(\widetilde T^m(x,u_0)).
\end{equation}

 Particularly, we take $u_0=c_2$, for every $M\in\mathbb{N}$ one has
	\begin{align*}
		S(f,\Lambda_{M}^{\vec{v}}(c_1,c_2))(x)=\sum_{m=0}^{M-1}\widetilde{f}(\widetilde{T}^m(x,c_2))
	\end{align*}
	for every $x\in X$. Then for $\mu$-a.e. $x\in X$,
\begin{align*}		
&\lim_{M\to\infty}\mathbb{E}(f,\Lambda_M^{\vec{v}}(c_1,c_2))(x)\\
=&\lim_{M\to\infty}\frac{1}{\#(\Lambda_{M}^{\vec{v}}(c_1,c_2))}S(f,\Lambda_{M}^{\vec{v}}(c_1,c_2))(x)\\
=&\lim_{M\to\infty}\frac{M}{\#(\Lambda_{M}^{\vec{v}}(c_1,c_2))}\frac{1}{M}\sum_{m=0}^{M-1}\widetilde{f}(\widetilde{T}^m(x,c_2))
=\frac{1}{c_2-c_1}\widehat f^{*}(x,c_2).
	\end{align*}
Set $f^*_{c_1,c_2}(x)=\frac{1}{c_2-c_1}\widehat f^*(x,c_2)$, then $\int_X f^*_{c_1,c_2}\,d\mu=\int_X f\,d\mu$.
	 \end{proof}

\begin{rem}\label{depend on b}
	Let $(X,\mathcal{B}_X,\mu,T)$ be a $\mathbb{Z}^2$-MPS, $\vec{v}=(1,\b)\in\mathbb{R}^2$. For every $f\in L^1(\mu)$, and $b>0$ we denote by $f^*_b=f_{-b/2,b/2}^*$. For a given $f\in L^1(\mu)$, there is an example showing that $f^*_b$ may depend on $b$ (see Example \ref{ex2}).
\end{rem}

By applying Theorem \ref{thm-1}, we can establish a directional version of the von Neumann theorem.
	\begin{thm}[Directional von Neumann ergodic theorem]\label{D-von Ne-thm} Let $(X,\mathcal{B}_X,\mu,T)$ be a $\mathbb{Z}^2$-MPS, $\vec{v}=(1,\beta)\in\mathbb{R}^2$ be a direction vector with $\beta \not\in \mathbb{Q}$. Then for every $c_1,c_2\in\mathbb{R}$ with $c_1<c_2$, $1\leq p<\infty$ and $f\in L^p(\mu)$ there exists $f^*_{c_1,c_2}\in L^p(\mu)$ such that
	$$\Vert \mathbb{E}(f,\Lambda_{M}^{\vec{v}}(c_1,c_2))-f^*_{c_1,c_2}\Vert_{p}\to 0, \ \text{as}\ \ M\to\infty.$$	
   	\end{thm}
   \begin{proof} For $f\in L^p(\mu)$, we only need to show that
   	$\{\mathbb{E}(f,\Lambda_{M}^{\vec{v}}(c_1,c_2))(x)\}_{M\in\mathbb{N}}$ is a Cauchy sequence in $L^p(\mu)$.
   	
   	For every $\varepsilon>0$, we take $g\in L^{\infty}(\mu)\subseteq L^1(\mu)$ such that $\Vert g-f\Vert_p<\varepsilon/4$. Then for any $M\in \mathbb{N}$,
  \begin{align*}
   		& \Vert\mathbb{E}(f,\Lambda_{M}^{\vec{v}}(c_1,c_2))-\mathbb{E}(g,\Lambda_{M}^{\vec{v}}(c_1,c_2))\Vert_p \\
   =& \bigg(\int\vert \mathbb{E}(f,\Lambda_M^{\vec{v}}(c_1,c_2))(x)-\mathbb{E}(g,\Lambda_{M}^{\vec{v}}(c_1,c_2))(x)\vert^pd\mu\bigg)^{1/p} \\
  \le & \bigg(\int_X\bigg|\frac{1}{\#\Lambda_M^{\vec{v}}(c_1,c_2)}\sum_{(m,n)\in\Lambda_M^{\vec{v}}(c_1,c_2)}\vert f(T^{(m,n)}x)-g(T^{(m,n)}x)\vert \bigg |^p d\mu(x)\bigg)^{1/p}\\
   \leq &\Vert f-g\Vert_p.
   	\end{align*}
By Theorem \ref{thm-1}, there exists $g^*\in L^1(\mu)$ such that
   	$$\mathbb{E}(g,\Lambda_{M}^{\vec{v}}(c_1,c_2))(x)\overset{\mu\text{-a.e.}}{\to} g^*(x), \ \ M\to\infty.$$
   	Clearly, $g^*\in L^{\infty}(\mu)$. 
   	Thus $g^*\in L^p(\mu)$ and $\big\vert \mathbb{E}(g,\Lambda_M^{\vec{v}}(c_1,c_2))(x)-g^*(x)\big\vert^p\to 0$ for $\mu$-a.e. $x\in X$. Hence,
   	$$\Vert \mathbb{E}(g,\Lambda_M^{\vec{v}}(c_1,c_2))-g^*\Vert_p \to 0,  \ \ M\to\infty.$$
   	Then there exists $M_0\in\mathbb{N}$, such that for every $M_1,M_2\geq M_0$ one has
   	$\Vert \mathbb{E}(g,\Lambda_{M_1}^{\vec{v}}(c_1,c_2))-\mathbb{E}(g,\Lambda_{M_2}^{\vec{v}}(c_1,c_2)) \Vert_p<\varepsilon/4$. Then
   	\begin{align*}
   		&\Vert\mathbb{E}(f,\Lambda_{M_1}^{\vec{v}}(c_1,c_2))-\mathbb{E}(f,\Lambda_{M_2}^{\vec{v}}(c_1,c_2)) \Vert_p\\
   \leq& \Vert\mathbb{E}(f,\Lambda_{M_1}^{\vec{v}}(c_1,c_2))-\mathbb{E}(g,\Lambda_{M_1}^{\vec{v}}(c_1,c_2))\Vert_p +
    \Vert\mathbb{E}(g,\Lambda_{M_1}^{\vec{v}}(c_1,c_2))-\mathbb{E}(g,\Lambda_{M_2}^{\vec{v}}(c_1,c_2))\Vert_p \\
   & +\Vert\mathbb{E}(g,\Lambda_{M_2}^{\vec{v}}(c_1,c_2))-\mathbb{E}(f,\Lambda_{M_2}^{\vec{v}}(c_1,c_2))\Vert_p \\
   \leq &2\Vert f-g\Vert_p+\Vert\mathbb{E}(g,\Lambda_{M_1}^{\vec{v}}(c_1,c_2))-\mathbb{E}(g,\Lambda_{M_2}^{\vec{v}}(c_1,c_2))\Vert_p<\varepsilon.
   	\end{align*}
   This ends our proof.
   \end{proof}
\subsection{Properties of directional ergodic average limits for different strips}

In this subsection, we will explore the behavior of the limit of directional ergodic averages when the strip changes.
\begin{prop}\label{continuous on b}
	Let $(X,\mathcal{B}_X,\mu,T)$ be a $\mathbb{Z}^2$-MPS and let
	$\vec{v}=(1,\beta)\in\mathbb{R}^2$ with $\beta\notin\mathbb{Q}$.
	For every $f\in L^1(\mu)$, the map
	$\Phi_f\colon (0,\infty)\to L^1(\mu)$ defined by $\Phi_f(b)=f_b^*$
	is continuous.
\end{prop}
\begin{proof}
	For $r>0$, write
	$a_M(r):=\#\Lambda_M^{\vec v}(r)$. It is easy to see that
	\begin{equation}\label{eq:strip-counting-continuity}
	\lim_{M\to\infty}\frac{a_M(r)}{M}=r.
	\end{equation}
	If $0<r<s$, then
	$\Lambda_M^{\vec v}(r)\subseteq\Lambda_M^{\vec v}(s)$. Hence, for $M\in \mathbb{N}$, using the $T$-invariance of $\mu$, we obtain
	\begin{align*}
	\norm{\mathbb E(f,\Lambda_M^{\vec v}(s))-\mathbb E(f,\Lambda_M^{\vec v}(r))}_1
	&\leq
	\frac{1}{a_M(s)}
	\sum_{(m,n)\in\Lambda_M^{\vec v}(s)\setminus\Lambda_M^{\vec v}(r)}
	\norm{f\circ T^{(m,n)}}_1 \\
	&\quad+
	\left|\frac1{a_M(s)}-\frac1{a_M(r)}\right|
	\sum_{(m,n)\in\Lambda_M^{\vec v}(r)}\norm{f\circ T^{(m,n)}}_1 \\
	&=2\frac{a_M(s)-a_M(r)}{a_M(s)}\norm{f}_1.
	\end{align*}
	Therefore, for arbitrary $b,b_0>0$,
	\begin{equation}\label{eq:finite-width-continuity-bound}
	\norm{\mathbb E(f,\Lambda_M^{\vec v}(b))-\mathbb E(f,\Lambda_M^{\vec v}(b_0))}_1
	\leq
	2\frac{|a_M(b)-a_M(b_0)|}
	{\max\{a_M(b),a_M(b_0)\}}\norm{f}_1
	\end{equation}
	for $M\in \mathbb{N}$.

	By Theorem~\ref{thm-1}, $\mathbb E(f,\Lambda_M^{\vec v}(b))\to f_b^*$ and
	$\mathbb E(f,\Lambda_M^{\vec v}(b_0))\to f_{b_0}^*$ almost everywhere. Thus Fatou's lemma,
	\eqref{eq:finite-width-continuity-bound}, and
	\eqref{eq:strip-counting-continuity} yield
	\begin{align*}
	\norm{f_b^*-f_{b_0}^*}_1
	&\leq \liminf_{M\to\infty}
	\norm{\mathbb E(f,\Lambda_M^{\vec v}(b))-\mathbb E(f,\Lambda_M^{\vec v}(b_0))}_1\\
	&\leq
	2\frac{|b-b_0|}{\max\{b,b_0\}}\norm{f}_1.
	\end{align*}
	The right-hand side tends to $0$ as $b\to b_0$. Hence
	$\Phi_f(b)=f_b^*$ is continuous from $(0,\infty)$ to $L^1(\mu)$.
\end{proof}
		

    \begin{prop}\label{continuous-on-f}
    Let $(X,\mathcal{B}_X,\mu,T)$ be a $\mathbb{Z}^2$-MPS, $\vec{v}=(1,\b)\in\mathbb{R}^2$ with $\b\notin\mathbb{Q}$. For every $b>0$, the map $\Phi_b\colon L^1(\mu)\to L^1(\mu)$ defined as $\Phi_b(f)=f^*_b$ is continuous.
    \end{prop}
    \begin{proof}
    	For every $g,f\in L^1(\mu)$, by Theorem \ref{thm-1}, there exist $g^*_b$ and $f^*_b\in L^1(\mu)$ satisfying (\ref{Brik-eq}) corresponding $g,f$. Then by Fatou's Lemma,
    	\begin{align}\label{e-continuous-f}
    		\Vert g^*_b-f^*_b\Vert_1&=\int_X\vert g_b^*-f_b^*\vert d\mu(x)\nonumber\\
    		&=\int_{X}\lim_{M\to\infty}\big\vert\mathbb{E}(g,\Lambda_M^{\vec{v}}(b))
    -\mathbb{E}(f,\Lambda_M^{\vec{v}}(b))\big\vert d\mu(x)\nonumber\\
    		&\leq \liminf_{M\to\infty}\int_X\frac{1}{\#\Lambda_M^{\vec{v}}(b)}\sum_{(m,n)\in\Lambda_M^{\vec{v}}(b)}\vert g(T^{(m,n)}x)-f(T^{(m,n)}x)\vert d\mu(x)\nonumber\\
    		&=\Vert g-f\Vert_1.
    	\end{align}
    	This establishes the continuity of $\Phi_b$.
    \end{proof}


  \begin{prop}\label{b-c}
	Let $(X,\mathcal{B}_X,\mu,T)$ be a $\mathbb{Z}^2$-MPS and let
	$\vec{v}=(1,\beta)\in\mathbb{R}^2$ with $\beta\notin\mathbb{Q}$.
	Then the following are equivalent:
\begin{itemize}
	\item[(1)] For every $b>0$ and every $f\in L^1(\mu)$,
	\[
	\mathbb{E}(f,\Lambda_M^{\vec{v}}(b))(x)\longrightarrow\int_X f\,d\mu
	\quad\text{as }M\to\infty
	\]
	for $\mu$-a.e. $x\in X$;
	\item[(2)] For every $c_1,c_2\in\mathbb{R}$ with $c_1<c_2$ and every
	$f\in L^1(\mu)$,
	\[
	\mathbb{E}(f,\Lambda_M^{\vec{v}}(c_1,c_2))(x)\longrightarrow\int_X f\,d\mu
	\quad\text{as }M\to\infty
	\]
	for $\mu$-a.e. $x\in X$.
\end{itemize}
\end{prop}
\begin{proof}
The implication $(2)\Rightarrow(1)$ follows by taking
$c_1=-b/2$ and $c_2=b/2$.

We prove $(1)\Rightarrow(2)$.
By decomposing $f$ into its positive and negative parts, it is enough to
consider $f\geq0$. Put
\[
b_0=c_2-c_1,
\qquad c_0=\frac{c_1+c_2}{2}=n_0+t_0,
\]
where $n_0\in\mathbb Z$ and $t_0\in[0,1)$.

If $t_0=0$, then for every $M$,
\[
\Lambda_M^{\vec v}(c_1,c_2)
 =\Lambda_M^{\vec v}(b_0)+(0,n_0).
\]
Hence, by the $T$-invariance of $\mu$ and (1),
\[
\mathbb E(f,\Lambda_M^{\vec v}(c_1,c_2))(x)
 =\mathbb E(f,\Lambda_M^{\vec v}(b_0))(T^{(0,n_0)}x)
 \longrightarrow \int_X f\,d\mu
\]
for $\mu$-a.e. $x$. Thus we may assume $t_0\in(0,1)$.

Choose a sequence $(\eta_i)_{i\ge1}$ such that
\[
0<\eta_i<1-t_0,\qquad \sum_{i=1}^{\infty}\eta_i<\infty,
\qquad \eta_i\longrightarrow0.
\]
Since $\beta$ is irrational, the set
$\{m\beta\pmod 1:m\in\mathbb N\}$ is dense in $\mathbb T$.
Thus, for every $i$, there exist $m_i\in\mathbb N, n_i\in\mathbb Z$ and $t_i\in (t_0, t_0+\eta_i)$
such that $n_i= m_i\beta+t_i$.
Set $\delta_i=t_i-t_0\in(0,\eta_i)$.
For $M>m_i$, define
\[
B_{M,i}=\Lambda_{M-m_i}^{\vec v}(b_0)+(m_i,n_i+n_0).
\]
If $(k,\ell)\in B_{M,i}$, then $m_i\leq k\leq M-1$ and
\[
k\beta+c_1+\delta_i
 \leq \ell\leq
k\beta+c_2+\delta_i.
\]
Consequently,
\begin{equation}\label{eq:b-c-symdiff}
\Lambda_M^{\vec v}(c_1,c_2)\triangle B_{M,i}
\subseteq L_i\cup \Lambda_M^{\vec v}(c_1,c_1+\eta_i)\cup \Lambda_M^{\vec v}(c_2,c_2+\eta_i),
\end{equation}
where
\[
L_i=\{(k,\ell)\in\Lambda^{\vec v}(c_1,c_2):0\leq k<m_i\}.
\]

By Theorem~\ref{thm-1}, for each $i$ there are nonnegative functions
$f^*_{c_1,c_1+\eta_i},f^*_{c_2,c_2+\eta_i}\in L^1(\mu)$ such that, for $\mu$-a.e. $x$,
\[
\mathbb E(f,\Lambda_M^{\vec v}(c_1,c_1+\eta_i))(x)\longrightarrow f^*_{c_1,c_1+\eta_i}(x),
\qquad
\mathbb E(f,\Lambda_M^{\vec v}(c_2,c_2+\eta_i))(x)\longrightarrow f^*_{c_2,c_2+\eta_i}(x),
\]
and
\begin{equation}\label{eq:b-c-boundary-integrals}
\int_X f^*_{c_1,c_1+\eta_i}\,d\mu=\int_X f^*_{c_2,c_2+\eta_i}\,d\mu=\int_X f\,d\mu.
\end{equation}

Since $\beta$ is irrational, it is easy to see that
\begin{equation}\label{eq:b-c-counting}
\frac{\#\Lambda_M^{\vec v}(c_1,c_2)}{M}\longrightarrow b_0,
\qquad
\frac{\#\Lambda_M^{\vec v}(c_1,c_1+\eta_i)}{M}\longrightarrow\eta_i,
\qquad
\frac{\#\Lambda_M^{\vec v}(c_2,c_2+\eta_i)}{M}\longrightarrow\eta_i,
\end{equation}
for every fixed $i$, and also
\[
\frac{\#\Lambda_{M-m_i}^{\vec v}(b_0)}{M}\longrightarrow b_0.
\]

Define
\[
Y_0=\left\{x\in X:
\lim_{M\to\infty}\mathbb E(f,\Lambda_M^{\vec v}(b_0))(x)
=\int_X f\,d\mu\right\}.
\]
For each $i\geq1$ and $j\in\{1,2\}$, put
\[
Y_{i,j}=\left\{x\in X:
\lim_{M\to\infty}\mathbb E
\bigl(f,\Lambda_M^{\vec v}(c_j,c_j+\eta_i)\bigr)(x)
=f^*_{c_j,c_j+\eta_i}(x)<\infty\right\},
\]
and let
\[
Z_i=\bigcap_{(m,n)\in L_i}(T^{(m,n)})^{-1}
\bigl(\{x\in X:|f(x)|<\infty\}\bigr).
\]
By assumption~(1), Theorem~\ref{thm-1}, and $f\in L^1(\mu)$,
\[
\mu(Y_0)=\mu(Y_{i,1})=\mu(Y_{i,2})=1.
\]
Since $L_i$ is finite and $\mu$ is $T$-invariant, $\mu(Z_i)=1$.
Consequently, the measurable set
\[
X_0=
\bigcap_{i=1}^{\infty}
\bigl((T^{(m_i,n_i+n_0)})^{-1}Y_0\cap Y_{i,1}\cap Y_{i,2}\cap Z_i\bigr)
\]
satisfies $\mu(X_0)=1$. For every $x\in X_0$ and every $i\geq1$,
all the limits defining $Y_0$ and $Y_{i,j}$ hold at the required
points, and $S(f,L_i)(x)<\infty$.

For such $x$ and fixed $i$, assumption (1) gives
\begin{align}
\frac{1}{M}S(f,B_{M,i})(x)
&=\frac{\#\Lambda_{M-m_i}^{\vec v}(b_0)}{M}
  \mathbb E(f,\Lambda_{M-m_i}^{\vec v}(b_0))(T^{(m_i,n_i+n_0)}x)
\longrightarrow b_0\int_X f\,d\mu.                                      \label{eq:b-c-main-translate}
\end{align}
Moreover, by \eqref{eq:b-c-counting},
\[
\frac1M S(f,\Lambda_M^{\vec v}(c_1,c_1+\eta_i))(x)\rightarrow\eta_i f^*_{c_1,c_1+\eta_i}(x),
\frac1M S(f,\Lambda_M^{\vec v}(c_2,c_2+\eta_i))(x)\rightarrow\eta_i f^*_{c_2,c_2+\eta_i}(x),
\]
while
\[
\frac1M S(f,L_i)(x)\rightarrow0.
\]
Since $f\ge0$, \eqref{eq:b-c-symdiff} and
\eqref{eq:b-c-main-translate} imply
\begin{equation}\label{eq:b-c-limsup}
\limsup_{M\to\infty}
\left|
\frac1M S(f,\Lambda_M^{\vec v}(c_1,c_2))(x)-b_0\int_X f\,d\mu
\right|
\leq
\eta_i\big(f^*_{c_1,c_1+\eta_i}(x)+f^*_{c_2,c_2+\eta_i}(x)\big)
\end{equation}
for every $i$.

Finally, by \eqref{eq:b-c-boundary-integrals} and Tonelli's theorem,
\[
\int_X\sum_{i=1}^{\infty}
\eta_i\big(f^*_{c_1,c_1+\eta_i}(x)+f^*_{c_2,c_2+\eta_i}(x)\big)\,d\mu
=2\int_X f\,d\mu\sum_{i=1}^{\infty}\eta_i<\infty.
\]
Hence
\[
\eta_i\big(f^*_{c_1,c_1+\eta_i}(x)+f^*_{c_2,c_2+\eta_i}(x)\big)\longrightarrow0
\]
for $\mu$-a.e. $x$. Let $i\to\infty$ in
\eqref{eq:b-c-limsup} yields
\[
\frac1M S(f,\Lambda_M^{\vec v}(c_1,c_2))(x)\longrightarrow b_0\int_X f\,d\mu
\]
for $\mu$-a.e. $x$. Using again
$\#\Lambda_M^{\vec v}(c_1,c_2)/M\to b_0$, we conclude that
\[
\mathbb E(f,\Lambda_M^{\vec v}(c_1,c_2))(x)\longrightarrow
\int_X f\,d\mu
\]
for $\mu$-a.e. $x$. This proves $(2)$ for nonnegative $f$, and the general
case follows by applying the result to $f^+$ and $f^-$.
\end{proof}

\section{Directional ergodicity and weak mixing}\label{sec: Der-Dwm}
In this section, based on the ergodic theorem from the previous section, we may introduce directional ergodicity and weak mixing. Then we will investigate the connections between the ergodic or weakly mixing properties of the original system, the directional action, and the corresponding skew-product action (see Theorem \ref{realtion-erg}, Theorem \ref{relation-wm} for a summary of the results).
\subsection{Directional ergodicity} This subsection is devoted to the introduction of directional ergodicity, the study of its basic properties, and the exploration of its relation to existing notions. Prior to this, we first present a proposition to clarify why invariant sets are not used to define directional ergodicity.

\begin{prop}\label{prop:reason}
Let $(X,\mathcal{B}_X,\mu,T)$ be a $\mathbb{Z}^2$-MPS and let
$\vec v=(1,\beta)$, where $\beta\notin\mathbb{Q}$.
For any $B\in\mathcal{B}_X$, the following are equivalent:
\begin{itemize}
    \item[(1)] $B$ is $T$-invariant modulo $\mu$, that is,
    $\mu(T^{-(m,n)}B\triangle B)=0$ for every $(m,n) \in\mathbb{Z}^2$;
    \item[(2)] for every $b>0$ and every $(m,n)\in\Lambda^{\vec v}(b)$,
    one has $\mu(T^{-(m,n)}B\triangle B)=0$;
    \item[(3)] there exists $b>0$ such that
    $\mu(T^{-(m,n)}B\triangle B)=0$ for every $(m,n)\in\Lambda^{\vec v}(b)$.
\end{itemize}
\end{prop}
\begin{proof}
The implications $(1)\Rightarrow(2)\Rightarrow(3)$ are immediate.
To prove $(3)\Rightarrow(1)$, fix $b>0$ and consider the additive
homomorphism
\[
    \ell:\mathbb{Z}^2\longrightarrow\mathbb{R},\qquad
    \ell(m,n)=n-\beta m.
\]
As $\beta$ is irrational, $\ell$ is injective and its image
$D=\mathbb{Z}+\beta\mathbb{Z}$ is dense in $\mathbb{R}$.
Moreover,
\[
    \Lambda^{\vec v}(b)
    =\{(m,n)\in\mathbb{Z}^2:|\ell(m,n)|\le b/2\}.
\]
We claim that every element of $\mathbb{Z}^2$ is a finite sum of
elements of $\Lambda^{\vec v}(b)$. Indeed, let $(m,n)\in\mathbb{Z}^2$,
without loss of generality put $t=\ell(m,n)>0$.
Choose an integer $N\ge2$ such that
$t/N<b/2$, and set $\eta=b/2-t/N>0$.
By density of $D$, we may choose $t_1,\ldots,t_{N-1}\in D$ satisfying
\[
    |t_j-t/N|<\eta/N\qquad (1\le j\le N-1).
\]
Set $t_N=t-\sum_{j=1}^{N-1}t_j\in D$. Then
\[
    |t_N-t/N|
    \le\sum_{j=1}^{N-1}|t_j-t/N|<\eta,
\]
and consequently $|t_j|<b/2$ for every $1\le j\le N$.
For each $j$, select $(m_j,n_j)\in\mathbb{Z}^2$ with
$\ell(m_j,n_j)=t_j$. Thus $(m_j,n_j)\in\Lambda^{\vec v}(b)$ and
\[
    \ell\Bigl((m,n)-\sum_{j=1}^N (m_j,n_j)\Bigr)
    =t-\sum_{j=1}^N t_j=0.
\]
Injectivity of $\ell$ yields $(m,n)=\sum_{j=1}^N (m_j,n_j)$, proving the claim.

Now let
\[
    H_B=\{(m,n)\in\mathbb{Z}^2:
       \mu(T^{-(m,n)}B\triangle B)=0\}.
\]
Since the action $T$ preserves $\mu$, $H_B$ is a subgroup of
$\mathbb{Z}^2$. Under assumption $(3)$ it contains
$\Lambda^{\vec v}(b)$, which generates $\mathbb{Z}^2$ by the
preceding argument. Hence $H_B=\mathbb{Z}^2$, which is $(1)$.
\end{proof}

With Birkhoff and von Neumann type ergodic theorems established, we define ergodicity along directions.

\begin{defn}\label{d-ergodicity}
	Let $(X,\mathcal{B}_X,\mu,T)$ be a $\mathbb{Z}^2$-MPS,
$\vec{v}=(1,\beta)\in\mathbb{R}^2$ be a direction vector with $\beta \not\in \mathbb{Q}$. We say $(X,\mathcal{B}_X,\mu,T)$ is
\textbf{$\vec{v}$-ergodic} if for every $b>0$ and $f\in L^1(\mu)$ one has
\begin{align}\label{def-v-ergodic}
	\mathbb{E}(f,\Lambda_M^{\vec{v}}(b))(x)\to \int_X fd\mu, \ \text{as}\ 	M\to \infty, \ \text{for}\ \mu\text{-a.e.}\ x\in X.
\end{align}
\end{defn}

Analogous to the general $\mathbb{Z}$-actions, we can characterize directional ergodicity in terms of subsets.

\begin{lem}\label{v-ergodic-set}
Let $(X,\mathcal{B}_X,\mu,T)$ be a $\mathbb{Z}^2$-MPS,
 $\vec{v}=(1,\beta)\in\mathbb{R}^2$ be a direction with $\beta \not\in \mathbb{Q}$.
 Let $\mathcal{D}$ be a semi-algebra that generates $\mathcal{B}_X$.
Then the following statements are equivalent:
\begin{enumerate}
\item[(1)] $(X,\mathcal{B}_X,\mu,T)$ is $\vec{v}$-ergodic;
\item[(2)] For any $A\in \mathcal{B}_X$, and $b>0$ one has
	$\mathbb{E}(1_A,\Lambda^{\vec{v}}_M(b))(x)\overset{M\to\infty}{\to}\mu(A)\ \text{for}\ \mu\text{-a.e}\ x\in X$;
\item[(3)] For any $A\in \mathcal{D}$, and $b>0$ one has $\mathbb{E}(1_A,\Lambda^{\vec{v}}_M(b))(x)\overset{M\to\infty}{\to}\mu(A)\ \text{for}\ \mu\text{-a.e}\ x\in X$;
\item[(4)] For any $A,B\in \mathcal{B}_X$, and $b>0$ one has
$$\lim_{M\to \infty}\frac{1}{\#(\Lambda_M^{\vec{v}}(b))}\sum_{(m,n)\in\Lambda_M^{\vec{v}}(b)} \mu (T^{-(m,n)}A\cap B)
=\mu(A)\mu(B);$$
\item[(5)] For any $A,B\in \mathcal{D}$, and $b>0$ one has 
$$\lim_{M\to \infty}\frac{1}{\#(\Lambda_M^{\vec{v}}(b))}\sum_{(m,n)\in\Lambda_M^{\vec{v}}(b)} \mu (T^{-(m,n)}A\cap B)
=\mu(A)\mu(B).$$
\end{enumerate}
\end{lem}
\begin{proof}
	$(1)\Rightarrow(2)\Rightarrow(3)$ and $(4)\Rightarrow(5)$ are trivial. We shall show $(2)\Rightarrow (1)$, $(2)\Leftrightarrow(4)$, and $(3)\Rightarrow(5)\Rightarrow (4)$.
	
	(2)$\Rightarrow$(1) For every $f\in L^1(\mu)$ and simple function $g\in L^1(\mu)$, by Theorem \ref{thm-1}, there exists $f^*_b, g^*_b\in L^1(\mu)$ satisfying
$$\lim_{M\to\infty}\mathbb{E}(f,\Lambda_M^{\vec{v}}(b))=f^*_{b},
\lim_{M\to\infty}\mathbb{E}(g,\Lambda_M^{\vec{v}}(b))=g^*_{b}.$$
It is straightforward to verify that (\ref{def-v-ergodic}) holds for $g$ by $(2)$, that is, $g^*_b=\int_X gd\mu$. Then by Proposition \ref{continuous-on-f},
	\begin{align*}
		&\Vert f^*_b-\int_Xfd\mu\Vert_1\leq \Vert f^*_b-g^*_b\Vert_1+\Vert \int_Xfd\mu-\int_Xgd\mu\Vert_1
		\overset{(\ref{e-continuous-f})}{\leq} 2 \cdot \Vert f-g\Vert_1.
	\end{align*}
	The density of simple functions in $L^1(\mu)$ yields $f^*_b=\int_Xfd\mu$	.
		
	(2)$\Rightarrow$(4) For every $A,B\in \mathcal{B}_X$, by $(2)$ and Dominated Convergence Theorem, one has
	\begin{align*}
		&\mu(A)\mu(B)\overset{(2)}{=}\int_B\lim_{M\to\infty}\mathbb{E}(1_A,\Lambda_M^{\vec{v}}(b))(x)d\mu =\lim_{M\to\infty}\int_B\mathbb{E}(1_A,\Lambda_M^{\vec{v}}(b))(x)d\mu \\
		=&\lim_{M\to\infty}\frac{1}{\#\Lambda_M^{\vec{v}}(b)}\sum_{(m,n)\in\Lambda_M^{\vec{v}}(b)}\int_B1_{T^{-(m,n)}A}d\mu
		=\lim_{M\to\infty}\frac{1}{\#\Lambda_M^{\vec{v}}(b)}\sum_{(m,n)\in\Lambda_M^{\vec{v}}(b)}\mu(T^{-(m,n)}A\cap B).
	\end{align*}
	
	(4)$\Rightarrow$(2) For any $A\in \mathcal{B}_X$, there exists $(1_A)_b^*\in L^1(\mu)$ satisfying
$$\lim_{M\to\infty}\mathbb{E}(1_A,\Lambda_M^{\vec{v}}(b))(x)=(1_A)_b^*(x).$$
Then for every $B\in\mathcal{B}_X$, by (4) and Dominated Convergence Theorem
	\begin{align*}
	&\int_B (1_A)_b^*(x) d\mu=\int_B\lim_{M\to\infty}\mathbb{E}(1_A,\Lambda_M^{\vec{v}}(b))(x)d\mu\\
	=&\lim_{M\to \infty}\frac{1}{\#(\Lambda_M^{\vec{v}}(b))}\sum_{(m,n)\in\Lambda_M^{\vec{v}}(b)} \mu (T^{-(m,n)}A\cap B) =\mu(A)\mu(B).
	\end{align*}
	The arbitrariness of $B$ implies that $(1_A)^*_b=\mu(A)$.
	
	(3)$\Rightarrow$(5) is similar to (2)$\Rightarrow$(4). (5)$\Rightarrow$(4) is similar to \cite[Theorem 1.17]{Peter}, we omit the proof.
\end{proof}

 For the relationship between directional ergodicity, ergodicity of the original $\mathbb{Z}^2$-action and that of the corresponding skew product, we have the following result.

 \begin{thm}\label{realtion-erg}
 		Let $(X,\mathcal{B}_X,\mu,T)$ be a $\mathbb{Z}^2$-MPS,
 $\vec{v}=(1,\beta)\in\mathbb{R}^2$ a direction vector with $\beta \not\in \mathbb{Q}$,
 and $(\widetilde{X},\widetilde{\mathcal{B}},\widetilde{\mu},\widetilde{T})$ the corresponding skew product. Then we have the following results:
 	\begin{itemize}
 		\item[(1)] If $(X,\mathcal{B}_X,\mu,T)$ is $\vec{v}$-ergodic, then $(X,\mathcal{B}_X,\mu,T)$ is ergodic.
 		\item [(2)] $(\widetilde{X},\widetilde{\mathcal{B}},\widetilde{\mu},\widetilde{T})$ is ergodic if and only if $(X,\mathcal{B}_X,\mu,T)$ is $\vec{v}$-ergodic.
 		\item [(3)] Ergodicity need not imply directional ergodicity.
 		\item [(4)] If $(X,\mathcal{B}_X,\mu,T)$ is ergodic, then there are at most countably many directions $\vec{v}$ such that it is not $\vec{v}$-ergodic.
 	\end{itemize}
 \end{thm}
 \begin{proof}(1) follows directly from the definition.
 	(2)  follows from Theorem \ref{cor:1}. (3)  follows from Example \ref{not-derg-Z^2-wm}. (4)  follows from Theorem \ref{thm:at most no ergodic}.
 \end{proof}

\begin{thm}\label{cor:1}
Let $(X,\mathcal{B}_X,\mu,T)$ be a $\mathbb{Z}^2$-MPS,
$\vec{v}=(1,\beta)\in\mathbb{R}^2$ a direction vector with $\beta \not\in \mathbb{Q}$. Let $(\widetilde{X},\widetilde{\mathcal{B}},\widetilde{\mu},\widetilde{T})$ be the corresponding skew product. Then $(\widetilde{X},\widetilde{\mathcal{B}},\widetilde{\mu},\widetilde{T})$ is ergodic if and only if $(X,\mathcal{B}_X,\mu,T)$ is $\vec{v}$-ergodic.
\end{thm}
\begin{proof}($\Leftarrow$) Suppose that $(X,\mathcal{B}_X,\mu,T)$ is $\vec{v}$-ergodic. Since $$\{B\times I\colon B\in\mathcal{B}_X, I\ \text{is a finite union of intervals}\}$$ forms a semi-algebra and generates the $\sigma$-algebra $\widetilde{\mathcal{B}}$. We only need to show that for every
$B\in\mathcal{B}_X$ and $a\in(0,1)$
\begin{align}\label{e-cor1-1}
	\lim_{M\to\infty}\frac{1}{M}\sum_{i=0}^{M-1}1_{B\times [0,a]}(\widetilde{T}^i(x,u))=a\cdot \mu(B)
\end{align}
for $\widetilde{\mu}$-a.e $(x,u)\in\widetilde{X}$. For this, we only need to show that for every $u_0\in [0,1)$ there exists a measurable $X_{u_0}\in\mathcal{B}_X$ with $\mu(X_{u_0})=1$ such that (\ref{e-cor1-1}) holds on $X_{u_0}\times\{u_0\}$.

We set $c_1=u_0-a$ and $c_2=u_0\in [0,1)$.
Since $(X,\mathcal{B}_X,\mu,T)$ is $\vec{v}$-ergodic, by Theorem \ref{thm-1} and Proposition \ref{b-c}, there exists a measurable $X_{u_0}\subseteq X$ with $\mu(X_{u_0})=1$ such that for every $x\in X_{u_0}$,
\begin{align}\label{e-cor1-2}
	\mathbb{E}(1_B,\Lambda_{M}^{\vec{v}}(c_1,c_2))(x)\to\mu(B),\ \text{as}\ M\to\infty.
\end{align}
We note that for every $M\in\mathbb{N}$ one has
\begin{align}\label{e-cor1-3}
	\sum_{i=0}^{M-1}1_{B\times[0,a]}(\widetilde{T}^i(x,u_0))=S(1_B,\Lambda_{M}^{\vec{v}}(c_1,c_2))(x).
\end{align}
Since $\beta\notin \mathbb{Q}$, we know that $\frac{\#\Lambda_M^{\vec{v}}(c_1,c_2)}{M}\to a$ as $M\to\infty$.
Combining with (\ref{e-cor1-2}) and (\ref{e-cor1-3}), we know that $(\ref{e-cor1-1})$ holds for $X_{u_0}\times \{u_0\}$.

($\Rightarrow$) Suppose $(\widetilde{X},\widetilde{\mathcal{B}},\widetilde{\mu},\widetilde{T})$ is ergodic. By Lemma \ref{v-ergodic-set}, it suffices to show that for every $A\in\mathcal{B}_X$ and $c_1,c_2\in \mathbb{R}$
with $c_1<c_2$,
\begin{align}\label{e-E-b}
	\lim_{M\to\infty}\mathbb{E}(1_A,\Lambda_M^{\vec{v}}(c_1,c_2))(x)=\mu(A)\ \text{for $\mu$-a.e.}\ x\in X.
\end{align}

By decomposing the strip $[c_1,c_2]$ into finitely many subintervals of length less than 1, and by the $T$-invariance of $\mu$, we can also assume
$0<c_2-c_1<1$ and $0\le c_2<1$.

We set $\omega=c_2-c_1$.
We define $\widetilde{f}\in L^1(\widetilde{\mu})$ as $\widetilde{f}(x,u)=1_A(x)1_{[0,\omega]}(u)$.
First assume  $u_0\in (0,1)$, for every $k\in\mathbb{N}$ we set
$0<t_k<\min\left\{\frac1{k},
	\sqrt{u_0},\sqrt{\omega},\sqrt{1-\omega}\right\}$ and
$$ \widetilde{g}_k(x,u)=1_A(T^{(0,1)}x)1_{[1-t_k^2,\ 1)}(u)+1_A(x)1_{[\omega-t_k^2,\ \omega]}(u)$$
for every $(x,u)\in\widetilde{X}$.
Applying the argument in Claim \ref{gk-f} of Theorem \ref{thm-1},
for every $k\in\mathbb{N}$ and $u\in(u_0-t_k^2,u_0), x\in X$, for any   $M\in\mathbb{N}$  one has
\begin{align}\label{qk-h}
\bigg\vert\sum_{i=0}^{M-1}\widetilde{f}(\widetilde{T}^i(x,u_0))
-\sum_{i=0}^{M-1}\widetilde{f}(\widetilde{T}^i(x,u))\bigg\vert
\leq\sum_{i=0}^{M-1}\widetilde{g}_k(\widetilde{T}^i(x,u)).
\end{align}

Then by the ergodicity of $(\widetilde{X},\widetilde{\mathcal{B}},\widetilde{\mu},\widetilde{T})$, there exist $T$-invariant $X_0\in \mathcal{B}_X$ with $\mu(X_0)=1$ and $u_k\in (u_0-t_k^2,u_0)$ for every $k\in\mathbb{N}$, such that for every $x\in X_0$
\begin{align}\label{e-qk-h-con}
	\lim_{M\to\infty}\frac{1}{M}\sum_{i=0}^{M-1}\widetilde{g}_k(\widetilde{T}^i(x,u_k))=2t_k^2\mu(A)\ \text{and}\ \lim_{M\to\infty}\frac{1}{M}\sum_{i=0}^{M-1}\widetilde{f}(\widetilde{T}^i(x,u_k))=\omega\cdot \mu(A)
\end{align}
Notice for every $x\in X_0$ and $k\in\mathbb{N}$ one has
\begin{align*}
&\limsup_{M\to\infty}\vert \frac{1}{ M}\sum_{i=0}^{M-1}\widetilde{f}(\widetilde{T}^i(x,u_0))-\omega\mu(A)\vert\\
\overset{(\ref{e-qk-h-con})}{=}
 &\limsup_{M\to\infty}\bigg\vert\frac{1}{M}\sum_{i=0}^{M-1}\widetilde{f}(\widetilde{T}^i(x,u_0))
-\frac{1}{ M}\sum_{i=0}^{M-1}\widetilde{f}(\widetilde{T}^i(x,u_k))\bigg\vert\\
	\overset{(\ref{qk-h})}{\leq}& \lim_{M\to\infty}\frac{1}{ M}\sum_{i=0}^{M-1}\widetilde{g}_k(\widetilde{T}^i(x,u_k))\overset{(\ref{e-qk-h-con})}{=}2t_k^2\mu(A).
\end{align*}
By the arbitrariness of $k$,
	$\lim_{M\to\infty}\frac{1}{ M}\sum_{i=0}^{M-1}\widetilde{f}(\widetilde{T}^i(x,u_0))=\omega\mu(A)$.

For $u_0=0$, $u_1=R_\beta(0)\in(0,1)$
and $\widetilde T(x,u_0)=(T^{(1,\lfloor\beta\rfloor)}x,u_1)$. By the preceding argument,
there exists a $T$-invariant $X_0\in \mathcal{B}_X$ with $\mu(X_0)=1$ such that for every $x\in X_0$
  $$\lim_{M\to\infty}\frac{1}{ M}\sum_{i=0}^{M-1}\widetilde{f}(\widetilde{T}^i(x,u_1))=\omega\mu(A)\ \text{as}\ M\to\infty.$$

For $\mu$-a.e. $x$ with $T^{(1,\lfloor\beta\rfloor)}x\in X_0$,
\[
\frac1M\sum_{i=0}^{M-1}\widetilde f(\widetilde T^i(x,0))
=
\frac{\widetilde f(x,0)}M
+
\frac{M-1}{M}
\left(
\frac1{M-1}\sum_{i=0}^{M-2}\widetilde f(\widetilde T^i(T^{(1,\lfloor\beta\rfloor)}x,u_1))
\right).
\]
Hence $\frac1M\sum_{i=0}^{M-1}\widetilde f(\widetilde T^i(x,0)) \rightarrow \omega\mu(A)$ as $M\to\infty.$

Then for any $u_0\in [0,1)$, by (\ref{eq:strip-skew-identity-thm1})
$$\lim_{M\to\infty}\mathbb{E}(1_A,\Lambda_M^{\vec{v}}(u_0-\omega,u_0))(x)=\lim_{M\to \infty} \frac{1}{ \omega M}\sum_{i=0}^{M-1}\widetilde{f}(\widetilde{T}^i(x,u_0))=\mu(A).$$
Particularly, we take $u_0=c_2$, we have
$\lim_{M\to\infty}\mathbb{E}(1_A,\Lambda_M^{\vec{v}}(c_1,c_2))(x)=\mu(A)$ for $c_2\in [0,1)$.

\end{proof}
	

\subsection{Directional weak mixing} Let $(X,\mathcal{B}_X,\mu,T)$ be a $\mathbb{Z}^2$-MPS, $\vec{v}=(1,\beta)\in\mathbb{R}^2$ be a direction vector with $\beta\not \in \mathbb{Q}$. We say $(X,\mathcal{B}_X,\mu,T)$ is \textbf{$\vec{v}$-weakly mixing} if the product $\mathbb{Z}^2$-action $(X\times X,\mathcal{B}_{X}\times \mathcal{B}_{X},\mu\times\mu,T\times T)$ is $\vec{v}$-ergodic. Clearly, $\vec{v}$-weakly mixing implies $\vec{v}$-ergodicity.

Analogous to the case of $\mathbb{Z}^2$-actions, the author of \cite{Liu} refers to a system as directionally weakly mixing if its directional Kronecker algebra is trivial. We will show that this definition is equivalent to ours. Specifically, for given direction $\vec{v}$ and $b>0$, the \textbf{directional Kronecker algebra} $\mathcal{K}^{\vec{v}}_{\mu}$ is defined as
 $$\mathcal{K}_{\mu}^{\vec{v}}\colon=\bigg\{B\in\mathcal{B}_X\colon\overline{\{U_{T^{(m,n)}}1_B\colon (m,n)\in \Lambda^{\vec{v}}(b)\}}\ \text{is compact in}\ L^2(\mu)\bigg\}.$$
For compact metrizable $\mathbb Z^2$-systems, the directional Kronecker
algebra is independent of the choice of $b>0$ and is a sub-$\sigma$-algebra;
see \cite{LX,Liu}. The next lemma extends these facts and the corresponding
Kronecker-algebra identity to arbitrary Lebesgue probability spaces.

\begin{lem}
\label{lem:compact-model-kron}
Let $(X,\mathcal B_X,\mu,T)$ be a $\mathbb Z^2$-MPS on a Lebesgue probability
space and let $\vec v=(1,\beta)$, where $\beta\notin\mathbb Q$. Then the family
$\mathcal K_\mu^{\vec v}$ defined above is a sub-$\sigma$-algebra modulo
$\mu$-null sets and is independent of $b>0$. Moreover, for the skew product
$(\widetilde X,\widetilde{\mathcal B},\widetilde\mu,\widetilde T)$ defined
in \eqref{skew-product-action}, one has
\begin{align}\label{Kron-d-Kron}
\mathcal K_{\widetilde\mu}
 =\mathcal K_\mu^{\vec v}\times\mathcal C
 \qquad\text{modulo }\widetilde\mu\text{-null sets},
\end{align}
where $\mathcal C$ denotes the Borel $\sigma$-algebra of $\mathbb T$.
\end{lem}
\begin{proof}
By Lemma \ref{0625}, there is a $\mathbb{Z}^2$-TDS $(Y,S)$ and a measure-theoretic isomorphism $\pi:(X,\mathcal{B}_X,\mu,T)\rightarrow(Y,\mathcal{B}_Y,\nu,S)$,
where $\nu$ is an $S$-invariant Borel probability measure and $\mathcal{B}_Y$ is the Borel $\sigma$-algebra of $Y$.

Since $\pi\circ T^{(m,n)}=S^{(m,n)}\circ \pi$ for any $(m,n)\in \mathbb{Z}^2$,
$\overline{\{U_S^g1_D:g\in\Lambda^{\vec v}(b)\}}$ is compact
if and only if $\overline{\{U_T^g1_{\pi^{-1}D}:g\in\Lambda^{\vec v}(b)\}}$ is compact for $D\in \mathcal{B}_Y$.
Since $\mathcal K_\nu^{\vec v}$ is independent of $b>0$,  $\mathcal K_\mu^{\vec v}$ is also independent of $b$
and  $\mathcal K_\mu^{\vec v} = \pi^{-1}\mathcal K_\nu^{\vec v}$ is a sub-$\sigma$-algebra.

Let
   \begin{align}
   	\widetilde{Y}\colon=Y\times \mathbb{T},\quad\widetilde{\nu}\colon=\nu\times m, \ \text{and}\ \ \widetilde{S}(y,u)\colon=(S^{(1,\lfloor\beta+u\rfloor)}y,\{\beta+u\})
   \end{align}
   for every $(y,u)\in\widetilde{Y}.$
By \cite[Lemma 4.4]{LX}, $\mathcal K_{\widetilde\nu}
 =\mathcal K_\nu^{\vec v}\times\mathcal C$.

The map $\pi\times\operatorname{id}$ is a measure-theoretic isomorphism
between $(\widetilde X,\widetilde\mu,
\widetilde T)$ and $(\widetilde{Y},\widetilde{\nu},\widetilde{S} )$:
\[
(\pi\times\operatorname{id})\circ\widetilde T
 =\widetilde S\circ(\pi\times\operatorname{id})
\qquad\widetilde\mu\text{-a.e.}
\]
Thus $$\mathcal K_{\widetilde\mu}=(\pi\times \operatorname{id})^{-1}(\mathcal K_{\widetilde\nu})
=(\pi\times \operatorname{id})^{-1}(\mathcal K_\nu^{\vec v}\times\mathcal C)
 =\mathcal K_\mu^{\vec v}\times\mathcal C.$$
\end{proof}
Using the same method in Lemma \ref{lem:compact-model-kron}, \cite[Theorem 1]{Liu} also applies
to \(\mathbb Z^2\)-MPS.
\begin{thm}\cite[Theorem 1]{Liu}\label{v-weakly mixing}
	Let  $(X,\mathcal{B}_X,\mu,T)$ be a $\mathbb{Z}^2$-MPS,
$\vec{v}=(1,\beta)\in\mathbb{R}^2$ be a direction vector with $\beta \not\in \mathbb{Q}$.
Then the following statements are equivalent:
\begin{itemize}
	\item [(1)]$\mathcal{K}_{\mu}^{\vec{v}}=\{\emptyset, X\}$;
	\item [(2)] For every $A,B\in\mathcal{B}_X$ and $b>0$, there exists $\mathcal{Q}_{A,B}\subset\Lambda^{\vec{v}}(b)$ with
        $$\lim_{M\to\infty}\frac{\#(\mathcal{Q}_{A,B}\cap\Lambda_M^{\vec{v}}(b) )}{\#\Lambda_M^{\vec{v}}(b)}=1$$
        such that for every sequence $\{(m_i,n_i)\}_{i\in\mathbb{N}}$ of $\mathcal{Q}_{A,B}$
        that satisfies $\{m_i\}_{i\in\mathbb{N}}$ is a strictly
 increasing sequence
        one has $$\lim_{i\to\infty}\mu(T^{-(m_i,n_i)}A\cap B)=\mu(A)\mu(B);$$
        \item[(3)] for every $A,B\in\mathcal{B}_X$ and $b>0$, one has
        $$\lim_{M\to\infty}\frac{1}{\#\Lambda_M^{\vec{v}}(b)}\sum_{(m,n)\in\Lambda_M^{\vec{v}}(b)}\big\vert\mu(T^{-(m,n)}A\cap B)-\mu(A)\mu(B)\big\vert=0.$$
\end{itemize}
\end{thm}

We now show that our notion of directional weak mixing is equivalent to the one defined by the first author in \cite{Liu} via the directional Kronecker algebra.

\begin{thm}\label{thm:weak mixing and ergodicity product}
	Let  $(X,\mathcal{B}_X,\mu,T)$ be a $\mathbb{Z}^2$-MPS,
$\vec{v}=(1,\beta)\in\mathbb{R}^2$ be a direction vector with $\beta \not\in \mathbb{Q}$. Then the following statements are equivalent.
	\begin{itemize}
		\item[(1)] $\mathcal{K}_{\mu}^{\vec{v}}=\{\emptyset, X\}$;
		\item[(2)] $(X,\mathcal{B}_X,\mu,T)$ is $\vec{v}$-weakly mixing;
		\item[(3)]  $(X\times X,\mathcal{B}_{X}\times \mathcal{B}_{X},\mu\times\mu,T\times T)$ is $\vec{v}$-weakly mixing.
	\end{itemize}
\end{thm}	
\begin{proof}

$(1) \Rightarrow (3)$ For every $A,B,C,D\in\mathcal{B}_X$, by Theorem \ref{v-weakly mixing} (2) there exist subsets $\mathcal{Q}_{A,B}$, $\mathcal{Q}_{C,D}\subset \Lambda^{\vec{v}}(b)$ with
\[\lim_{M\to \infty}\frac{\#(\mathcal{Q}_{A,B}\cap\Lambda_M^{\vec{v}}(b))}{\#\Lambda_M^{\vec{v}}(b)}
=\lim_{M\to \infty}\frac{\#(\mathcal{Q}_{C,D}\cap\Lambda_M^{\vec{v}}(b))}{\#\Lambda_M^{\vec{v}}(b)}=1\]
such that for any subset $\{(m_i,n_i)\}_{i=1}^\infty$ of $\mathcal{Q}_{A,B}\cap \mathcal{Q}_{C,D}$
that satisfies $\{m_i\}_{i\in\mathbb{N}}$ is a strictly
monotone increasing sequence
\[\lim_{i\to\infty}\mu(T^{-(m_i,n_i)}A\cap B)=\mu(A)\mu(B),\]
\[\lim_{i\to\infty}\mu(T^{-(m_i,n_i)}C\cap D)=\mu(C)\mu(D).\]
Then $(\mu\times\mu)\big((T\times T)^{-(m_i,n_i)}(A\times C)\cap (B\times D)\big)\to(\mu\times\mu)(A\times C)(\mu\times\mu)(B\times D)$ as $i\to\infty$.
Note that $\lim_{M\to \infty}\frac{\#(\mathcal{Q}_{A,B}\cap\mathcal{Q}_{C,D}\cap\Lambda_M^{\vec{v}}(b))}{\#\Lambda_M^{\vec{v}}(b)}=1$.
Then it is easy to see that
\begin{align*}
& \lim_{M\to \infty}\frac{1}{\#(\Lambda_M^{\vec{v}}(b))}\sum_{(m,n)\in\Lambda_M^{\vec{v}}(b)}
|(\mu\times\mu)\big((T\times T)^{-(m,n)}(A\times C)\cap (B\times D)\big)\\
&-(\mu\times\mu)(A\times C)(\mu\times\mu)(B\times D)|=0.
\end{align*}

Since $\{A\times C: A,C\in \mathcal{B}_X\}$ forms a semi-algebra that generates $\mathcal{B}_{X}\times \mathcal{B}_{X}$,
then the above equality holds for all measurable sets in $\mathcal{B}_{X}\times \mathcal{B}_{X}$.
Then by Theorem \ref{v-weakly mixing} (3), $(X\times X,\mathcal{B}_{X}\times \mathcal{B}_{X},\mu\times\mu,T\times T)$ is $\vec{v}$-weakly mixing.

$(3)\Rightarrow (2)$ is trivial.

$(2)\Rightarrow(1)$ By Theorem \ref{v-weakly mixing} (3), it is sufficient to show that for every $A,B\in\mathcal{B}_X$, we have
\begin{align}\label{e-wm-pro-er}
	\lim_{M\to \infty}\frac{1}{\#(\Lambda_M^{\vec{v}}(b))}\sum_{(m,n)\in\Lambda_M^{\vec{v}}(b)} \vert\mu (T^{-(m,n)}A\cap B)-\mu(A)\mu(B)\vert^2=0.
\end{align}
Since $(X\times X,\mathcal{B}_{X}\times \mathcal{B}_{X},\mu\times\mu,T\times T)$ is $\vec{v}$-ergodic, by Lemma \ref{v-ergodic-set}, we have
\begin{align}\label{e-wm-pro-er1}
	\mu(A)\mu(B)&
	=\lim_{M\to \infty}\int\mathbb{E}(1_{A\times X},\Lambda_M^{\vec{v}}(b))(x,y)1_{B\times X}(x,y)d\mu\times \mu(x,y)\nonumber\\
	&=\lim_{M\to \infty}\frac{1}{\#\Lambda_M^{\vec{v}}(b)}\sum_{(m,n)\in\Lambda_M^{\vec{v}}(b)} \mu (T^{-(m,n)}A\cap B),
\end{align}
and
\begin{align}\label{e-wm-pro-er2}
  \mu(A)^2\mu(B)^2
  &=\lim_{M\to\infty}\int \mathbb{E}(1_{A\times A},\Lambda_M^{\vec{v}}(b))(x,y)1_{B\times B}(x,y)d\mu\times \mu(x,y) \nonumber\\
&=\lim_{M\to \infty}\frac{1}{\#\Lambda_M^{\vec{v}}(b)}\sum_{(m,n)\in\Lambda_M^{\vec{v}}(b)} \mu (T^{-(m,n)}A\cap B)^2.
\end{align}
Then (\ref{e-wm-pro-er}) follows from (\ref{e-wm-pro-er1}) and (\ref{e-wm-pro-er2}).

\end{proof}

\begin{rem} By the above analysis and Theorem \ref{v-weakly mixing}(2), it is not hard to prove that if $(X,\mathcal{B}_X,\mu,T)$ and $(Y,\mathcal{B}_Y,\nu,S)$ are $\vec{v}$-weakly mixing then
$(X\times Y,\mathcal{B}_{X}\times \mathcal{B}_{Y},\mu\times \nu,T\times S)$ is also $\vec{v}$-weakly mixing.
\end{rem}

Similar to directional ergodicity, we shall examine the relationship between directional weak mixing, weak mixing of the original $\mathbb{Z}^2$-action and that of the corresponding skew product. We summarize the results as follows.
\begin{thm}\label{relation-wm}
	Let $(X,\mathcal{B}_X,\mu,T)$ be a $\mathbb{Z}^2$-MPS,
$\vec{v}=(1,\beta)\in\mathbb{R}^2$ be a direction vector with $\beta \not\in \mathbb{Q}$, and $(\widetilde{X},\widetilde{\mathcal{B}},\widetilde{\mu},\widetilde{T})$ the corresponding skew product. Then we have:
	\begin{itemize}
	\item[(1)] $(X,\mathcal{B}_X,\mu,T)$ is $\vec{v}$-weakly mixing if and only if it is weakly mixing and $\vec{v}$-ergodic.
	\item[(2)] $(\widetilde{X},\widetilde{\mathcal{B}},\widetilde{\mu},\widetilde{T})$ is not weakly mixing.
	\item[(3)] Weak mixing need not imply directional weak mixing.
	\item[(4)] If $(X,\mathcal{B}_X,\mu,T)$ is weakly mixing, then there are at most countably many $\vec{v}$
such that it is not $\vec{v}$-weakly mixing.
	\end{itemize}	
\end{thm}
\begin{proof}
	(1) follows from Theorem \ref{thm:weak mixing+directional ergodic=directional weak}.
	(2) follows from Theorem \ref{X-weakly mixing} and (\ref{Kron-d-Kron}). (3) refers to Example \ref{not-derg-Z^2-wm}. (4) refers to Theorem \ref{thm:at most no ergodic}.
\end{proof}
We need the following lemma.
\begin{lem}\label{d-wm and int con}
Let $(X,\mathcal{B}_X,\mu,T)$ be a $\mathbb{Z}^2$-MPS,
$\vec{v}=(1,\beta)\in\mathbb{R}^2$ be a direction vector with $\beta \not\in \mathbb{Q}$. The following are equivalent:
	\begin{itemize}
		\item [(1)]$(X,\mathcal{B}_X,\mu,T)$ is $\vec{v}$-weakly mixing;
		\item [(2)] For every eigenfunction $\widetilde{f}$ of $(\widetilde{X},\widetilde{\mathcal{B}},\widetilde{\mu},\widetilde{T})$, there exists a constant $c_f$ such that $\int_{\mathbb{T}} \widetilde{f}(x,u)dm(u)=c_f$ for $\mu$-a.e. $x\in X$.
	\end{itemize}	
\end{lem}
\begin{proof}(1)$\Rightarrow$ (2).
 By (\ref{kro-eigen}), every eigenfunction $\widetilde{f}$ of $(\widetilde{X},\widetilde{\mathcal{B}},\widetilde{\mu},\widetilde{T})$ is $\mathcal{K}_{\widetilde{\mu}}$-measurable. By (\ref{Kron-d-Kron}), $\int_{\mathbb{T}}\widetilde{f}(x,u)dm(u)$ is $\mathcal{K}_{\mu}^{\vec{v}}$-measurable. Since $(X,\mathcal{B}_X,\mu,T)$ is $\vec{v}$-weakly mixing,
 by Theorem \ref{thm:weak mixing and ergodicity product} $\mathcal{K}_{\mu}^{\vec{v}}=\{\emptyset,X\}$. Hence $\int_{\mathbb{T}} \widetilde{f}(x,u)dm(u)$ is constant $\mu$-a.e.
	
	(2)$\Rightarrow$ (1). We shall show that $\mathcal{K}_{\mu}^{\vec{v}}$ is trivial. For $A\in\mathcal{K}_{\mu}^{\vec{v}}$, we set $\widetilde{g}(x,u)=1_A(x)1_{\mathbb{T}}(u)$ for every $x\in X$ and $u\in\mathbb{T}$. By (\ref{Kron-d-Kron}), $\widetilde{g}$ is $\mathcal{K}_{\widetilde{\mu}}$-measurable and $\int_{\mathbb{T}}\widetilde{g}(x,u)dm(u)=1_A(x)$ for $\mu$-a.e. $x\in X$. By (\ref{kro-eigen}), 
	there exist eigenfunctions $\widetilde{f}_{n,i}$ of $(\widetilde{X},\widetilde{\mathcal{B}},\widetilde{\mu},\widetilde{T})$
with $U_{\widetilde{T}}\widetilde{f}_{n,i}=\lambda_{n,i} \widetilde{f}_{n,i}$ for some $\lambda_{n,i}\in\mathbb{C}$
and $a_{n,i}\in \mathbb{C}$ such that
	$$\widetilde{h}_n\colon=\sum_{i=1}^na_{n,i}\widetilde{f}_{n,i}\to \widetilde{g}\ \ \text{in}\  L^2(\widetilde{\mu}),\ \text{as}\   n\to\infty.$$ By (2), there exist constants $c_{n,i}$ such that $\int_{\mathbb{T}} \widetilde{f}_{n,i}(x,u)dm(u)=c_{n,i}$ for $\mu$-a.e. $x\in X$ and every $n\in\mathbb{N}$. Note that
	\begin{align*}
		&\limsup_{n\to\infty}\int_{X}\bigg\vert\int_{\mathbb{T}}\widetilde{h}_n(x,u)dm(u)-\int_{\mathbb{T}}\widetilde{g}(x,u)dm(u)\bigg\vert d\mu(x)\\
		\leq &\limsup_{n\to\infty}\int_{X}\int_{\mathbb{T}}\vert \widetilde{h}_n(x,u)-\widetilde{g}(x,u)\vert dm(u)d\mu(x)\\
		\leq &\limsup_{n\to\infty}\bigg(\int_{\widetilde{X}}\vert \widetilde{h}_n(x,u)-\widetilde{g}(x,u)\vert^2d\widetilde{\mu}\bigg)^{1/2}=0.
	\end{align*}
	Then there exists a subsequence $\{\widetilde h_{n_j}\}_{j\in\mathbb{N}}$ such that
	$$\sum_{i=1}^{n_j}a_{n_j,i}c_{n_j,i}=\int_{\mathbb{T}}\widetilde{h}_{n_j}(x,u)dm(u)\to \int_{\mathbb{T}}\widetilde{g}(x,u)dm(u)=1_A(x)\ \text{as}\ j\to\infty,$$
	for $\mu$-a.e. $x\in X$. This implies $\mu(A)=0$ or 1. This ends our proof. 	
\end{proof}

\begin{thm}\label{thm:weak mixing+directional ergodic=directional weak}
	Let $(X,\mathcal{B}_X,\mu,T)$ be a $\mathbb{Z}^2$-MPS,
$\vec{v}=(1,\beta)\in\mathbb{R}^2$ be a direction vector with $\beta \not\in \mathbb{Q}$.
Then the following statements are equivalent:
	\begin{itemize}
		\item[(1)] $(X,\mathcal{B}_X,\mu,T)$ is weakly mixing and $\vec{v}$-ergodic;
		\item[(2)]  $(X,\mathcal{B}_X,\mu,T)$ is $\vec{v}$-weak mixing.
	\end{itemize}
\end{thm}
\begin{proof}
 $(2)\Rightarrow (1)$ By definition,  $\vec{v}$-weak mixing implies $\vec{v}$-ergodicity.
 Moreover, for every $b>0$,
	$$\overline{\{U_{T^{(m,n)}}1_B\colon (m,n)\in\Lambda^{\vec{v}}(b)\}}\subseteq \overline{\{U_{T^{(m,n)}}1_B\colon (m,n)\in\mathbb{Z}^2\}}.$$ Hence, by Theorem \ref{thm:weak mixing and ergodicity product} and Theorem \ref{X-weakly mixing}, $\vec{v}$-weak mixing implies weak mixing.

 $(1)\Rightarrow (2)$ By Lemma \ref{d-wm and int con}, it suffices to show that for every eigenfunction $\widetilde{f}$ of $(\widetilde{X},\widetilde{\mathcal{B}},\widetilde{\mu},\widetilde{T})$, there exists some constant $c_f$ such that $\int_{\mathbb{T}} \widetilde{f}(x,u)dm(u)=c_f$ for $\mu$-a.e. $x\in X$.

 Let $\widetilde{f}$ be an eigenfunction of $(\widetilde{X},\widetilde{\mathcal{B}},\widetilde{\mu},\widetilde{T})$ corresponding to eigenvalue $\lambda$. For every $(m,n)\in\mathbb{Z}^2$, we set
 $$\widetilde{g}_{m,n}(x,u)=\widetilde{f}(T^{(m,n)}x,u)$$
 for every $(x,u)\in\widetilde{X}$. We note that $\widetilde{T}\circ (T^{(m,n)}\times id)=(T^{(m,n)}\times id)\circ\widetilde{T}$ for every $(m,n)\in\mathbb{Z}^2$. Hence
 	\begin{align*}
 	U_{\widetilde T}\widetilde g_{m,n} =&
 	\widetilde f\circ (T^{(m,n)}\times \mathrm{id})\circ \widetilde T
 	=
 	\widetilde f\circ \widetilde T\circ (T^{(m,n)}\times \mathrm{id}) \\
 	=&
 	\lambda\,\widetilde f\circ (T^{(m,n)}\times \mathrm{id})
 	=
 	\lambda\,\widetilde g_{m,n}.
 \end{align*}
 This means $\widetilde{g}_{m,n}$ is also an eigenfunction of $(\widetilde{X},\widetilde{\mathcal{B}},\widetilde{\mu},\widetilde{T})$ corresponding to eigenvalue $\lambda$. Since $(X,\mathcal{B}_X,\mu,T)$ is $\vec{v}$-ergodic, $(\widetilde{X},\widetilde{\mathcal{B}},\widetilde{\mu},\widetilde{T})$ is ergodic by Theorem \ref{cor:1}. Then there exists a constant $c_{m,n}$ such that $\widetilde{g}_{m,n}=c_{m,n}\widetilde{f}$. In fact, $\vert c_{m,n}\vert=1$, because of $\Vert \widetilde{g}_{m,n}\Vert_2=\Vert\widetilde{f}\Vert_2$.

 We set $f(x)\colon=\int_{\mathbb{T}}\widetilde{f}(x,u)dm(u)$. A direct verification shows that $c\colon \mathbb{Z}^2 \to\{z\colon\vert z\vert=1\}, (m,n)\mapsto c_{m,n}$ is a homomorphism.
 If $f=0$  $\mu$-a.e., then $f$ is constant $\mu$-a.e. If $f$ is not identical to zero,
 then $f$ is an eigenfunction of $(X,\mathcal{B}_X,\mu,T)$.
By the weak mixing of $(X,\mathcal{B}_X,\mu,T)$ and Theorem \ref{X-weakly mixing}, $c_{m,n}=1$ and $f$ is constant $\mu$-a.e.
\end{proof}

\section{Directional ergodicity and weak mixing via spectrum}\label{sec:ergodicity and mixing via spectrum}
In this section, we investigate directional ergodicity and weak mixing via spectrum introduced by Robinson, Rosenblatt and \c{S}ahin \cite{RRS2023}. We prove that these spectral notions coincide with the corresponding notions defined in this paper. In particular, combined with our earlier result, we give an affirmative answer to the question posed in \cite[Section 7]{RRS2023}: whether the spectral weak mixing considered there is equivalent to Liu's directional weak mixing \cite{Liu}.

\subsection{Definition of ergodicity and weak mixing via spectrum}
In this subsection, we recall the definition of ergodicity and weak mixing introduced by Robinson, Rosenblatt and \c{S}ahin \cite{RRS2023} from the spectrum viewpoint.

Let $(X,\mathcal{B}_X,\mu,T)$ be a $\mathbb{Z}^2$-MPS. We can define its \textbf{$\mathbb{R}^2$-suspension flow $(\widetilde{X}',\widetilde{\mathcal{B}}',\widetilde{\mu}',$ $(\Phi_{t,s})_{t,s\in\mathbb{R}})$} as follows. Let $\widetilde{X}':=\widetilde{X}\times \mathbb{T}=X\times\mathbb{T}\times\mathbb{T}$ be the product space of $\widetilde{X}$ and $\mathbb{T}$, and $\widetilde{\mathcal{B}}'$ the $\sigma$-algebra generated by $\widetilde{\mathcal{B}}$ and $\mathcal{C}$. Let $\widetilde{\mu}':=\widetilde{\mu}\times m=\mu\times m\times m$, where $m$ is the Lebesgue measure on $\mathbb{T}$. For every $(t,s)\in\mathbb{R}^2$, we set
\[\Phi_{t,s}:\widetilde{X}'\to \widetilde{X}',~~(x,u,v)\mapsto \big(T^{(\lfloor s+v\rfloor,\lfloor t+u\rfloor)}x,\{u+t\},\{v+s\}\big).\]
It is straightforward to verify that $(\Phi_{t,s})_{t,s\in\mathbb{R}}$ is a $\widetilde{\mathcal{B}}'$-measurable $\mathbb{R}^2$-action preserving $\widetilde{\mu}'$. Note that for every $n\in\mathbb{Z}$, $\b\in\mathbb{R}$ and $(x,u,v)\in\widetilde{X}'$,
\begin{align}\label{flow-skew product action}
	\Phi_{n\b,n}(x,u,v)&=(T^{(\lfloor n+v\rfloor,\lfloor n\b+u \rfloor)}x,\{n\b +u\},\{n+v\})\nonumber\\
	&=\big(T^{(n,\lfloor n\b+u\rfloor)}x,\{n\b +u\},v\big)=\big(\widetilde{T}^n(x,u),v\big)
\end{align}
where $\widetilde{T}$ is the $\mathbb{Z}$-action on the skew product  corresponding to $(1,\b)$ (see (\ref{skew-product-action})).
Define a closed $\Phi$-invariant Hilbert space
\[L_{\mathcal{O}}^2(\widetilde{X}',\widetilde{\mu}'):=\mathcal{O}^{\perp}\subset L^2(\widetilde{X}',\widetilde{\mu}'),\]
where
\[\mathcal{O}:=\big\{F(x,u,v)=f(u,v):f\in L^2(\mathbb{T}^2,m\times m)\big\}\subset L^2(\widetilde{X}',\widetilde{\mu}').\]

\begin{defn}
Let $(X,\mathcal{B}_X,\mu,T)$ be a $\mathbb{Z}^2$-MPS, $(\widetilde{X}',\widetilde{\mathcal{B}}',\widetilde{\mu}',(\Phi_{t,s})_{t,s\in\mathbb{R}})$ be its $\mathbb{R}^2$ suspension flow and $\vec{v}=(1,\beta)$ be a direction vector with $\beta \not \in \mathbb{Q}$. We say $T$ is
\begin{enumerate}
    \item \textbf{$\vec{v}$-ergodic in terms of spectrum} if the $\mathbb{R}$-flow $(\widetilde{X}',\widetilde{\mathcal{B}}',\widetilde{\mu}',(\Phi_{s\b,s})_{s\in\mathbb{R}})$ has no non-zero invariant function in $L_{\mathcal{O}}^2(\widetilde{X}',\widetilde{\mu}')$;
    \item \textbf{$\vec{v}$-weakly mixing in terms of spectrum} if the $\mathbb{R}$-flow $(\widetilde{X}',\widetilde{\mathcal{B}}',\widetilde{\mu}',(\Phi_{s\b,s})_{s\in\mathbb{R}})$ has no non-zero eigenfunction in $L_{\mathcal{O}}^2(\widetilde{X}',\widetilde{\mu}')$.
\end{enumerate}
\end{defn}

\subsection{Directional ergodicity in terms of spectrum}
In this subsection, we prove that directional ergodicity is equivalent to directional ergodicity in terms of spectrum.
\begin{thm}\label{thm:ve=>ves}
   Let $(X,\mathcal{B}_X,\mu,T)$ be a $\mathbb{Z}^2$-MPS
   and $\vec{v}=(1,\beta)$ be a direction vector with $\beta \not\in \mathbb{Q}$.
   Then the following statements are equivalent:
   \begin{enumerate}
        \item[(1)] $(X,\mathcal{B}_X,\mu,T)$ is $\vec{v}$-ergodic;
		\item[(2)] $(X,\mathcal{B}_X,\mu,T)$ is $\vec{v}$-ergodic in terms of spectrum;
		\item[(3)] If $f\in L^2(\widetilde{X}',\widetilde{\mu}')$ and $f\circ\Phi_{s\beta,s}=f$ for any $s\in \mathbb{R}$, then $f$ is constant, $\widetilde{\mu}'$-a.e.
	\end{enumerate}
\end{thm}
\begin{proof}
Without loss of generality, we assume that $\beta>0$.

(1) $\Rightarrow$ (2) Let $f\in L_{\mathcal{O}}^2(\widetilde{X}',\widetilde{\mu}')$
be an invariant function with respect to the $\mathbb{R}$-flow $\{\Phi_{s\b,s}\}_{s\in\mathbb{R}}$. We shall show that $f$ is zero $\widetilde{\mu}'$-a.e.

For every $v\in\mathbb{T}$, we define $h_{v}(x,u)=f(x,u,v)$ for every $(x,u)\in\widetilde{X}$. Then for $m$-a.e $v\in\mathbb{T}$, one has $h_v\in L^1(\widetilde{\mu})$. For every $n\in\mathbb{Z}$,
\begin{align*}
	h_v(\widetilde{T}^n(x,u))=f(\widetilde{T}^n(x,u),v)\overset{(\ref{flow-skew product action})}{=}f(\Phi_{n\b,n}(x,u,v))=f(x,u,v)=h_v(x,u)
\end{align*}
for $\widetilde{\mu}$-a.e. $(x,u)\in\widetilde{X}$.

Since $(X,\mathcal{B}_X,\mu,T)$ is $\vec{v}$-ergodic, by Theorem \ref{cor:1},
the skew product $(\widetilde{X},\widetilde{\mathcal{B}},\widetilde{\mu},\widetilde{T})$ is ergodic. Hence $h_v$ is $\widetilde{\mu}$-a.e. constant, which implies $f\in \mathcal{O}$. Noting that $f\in L_{\mathcal{O}}^2(\widetilde{X}',\widetilde{\mu}')$, one has $f=0$ $\widetilde{\mu}'$-a.e.

(2) $\Rightarrow$ (3)
Let $f\in L^2(\widetilde{X}',\widetilde{\mu}')$ such that $f\circ\Phi_{s\b,s}=f$ for any $s\in \mathbb{R}$. We define $g\in L^2(m\times m)$ by $g(u,v)=\int_X f(x,u,v)d\mu(x)$ for $m\times m$-a.e. $(u,v)\in \mathbb{T}^2$. Since $\mu$ is $T$-invariant and $f$ is invariant under $\Phi_{s\b,s}$ for every $s\in\mathbb{R}$ one has
\begin{align*}
	g(u,v)=&\int_X f(x,u,v)d\mu(x)=\int f(\Phi_{s\b,s}(x,u,v))d\mu(x)\\
	=&\int_X f(T^{(\lfloor s+v\rfloor,\lfloor s\beta+u\rfloor)}x,\{u+s\b\},\{v+s\}) d\mu(x)\\
	=&\int_X f(x,\{u+s\b\},\{v+s\}) d\mu(x) =g(\{u+s\b\},\{v+s\}).
	\end{align*}
Since $\beta \in \mathbb{R}\setminus \mathbb{Q}$, $g(u,v)$ is $m\times m$-a.e. equal to a constant,
which we denote by $C$. Note that, for every $h(u,v)\in \mathcal{O}$,
\begin{align*}
	&\int_{\widetilde{X}'}h(u,v)(\overline{f(x,u,v)-C})d\widetilde{\mu}'\\
=&\int_{\mathbb{T}^2}h(u,v)\int_X(\overline{f(x,u,v)-C})d\mu(x)dm\times m(u,v)=0,
\end{align*}
which means $f-C\in L_{\mathcal{O}}^2(\widetilde{X}',\widetilde{\mu}')$. Moreover,
$(f-C)\circ\Phi_{s\b,s}=f-C$ for every $s\in\mathbb{R}$. Since $(X,\mathcal{B}_X,\mu,T)$ is $\vec{v}$-ergodic in terms of spectrum, $f-C=0$, i.e. $f(x,u,v)$ is $\widetilde{\mu}'$-a.e constant.

(3)$\Rightarrow$(1) By Theorem \ref{cor:1}, we only need to show that $(\widetilde{X},\widetilde{\mathcal{B}},\widetilde{\mu},\widetilde{T})$ is ergodic. Thus, it is sufficient to show that for every $B\in\mathcal{B}_X$ and $a\in(0,1)$ one has
\begin{align}\label{e-413}
	\lim_{N\to\infty}\frac{1}{N}\sum_{k=0}^{N-1}1_{B\times [0,a]}(\widetilde{T}^{k}(x,u))=a \cdot \mu(B)
\end{align}
for $\widetilde{\mu}$-a.e. $(x,u)\in \widetilde{X}$.

For every $n \in\mathbb{N}$ such that
$\frac{2\beta}{n}<\min\{a,1-a\}$, by Birkhoff's pointwise ergodic theorem and (3), we have
\begin{align}\label{e4-11}
	\frac{1}{N}\int _0^N1_{B\times [\frac{2\beta}{n},a]\times (\frac{1}{n},\frac{2}{n})}\big(\Phi_{s\beta,s}(x,u,v)\big)ds\to \frac{\mu(B)}{n}\cdot(a-\frac{2\beta}{n})\ \text{as}\ N\to\infty,
\end{align}
and
\begin{align}\label{e4-12}
	\frac{1}{N}\int _0^N1_{B\times [0,a+\frac{2\beta}{n}]\times (\frac{1}{n},\frac{2}{n})}\big(\Phi_{s\beta,s}(x,u,v)\big)ds\to \frac{\mu(B)}{n}\cdot(a+\frac{2\beta}{n})\ \text{as}\ N\to\infty,
\end{align}
for $\widetilde{\mu}'$-a.e. $(x,u,v)\in \widetilde{X}'$. Then there exists $v_n\in(0,1/n^2)$ and $\widetilde{X}_n\subseteq\widetilde{X}$ with $\widetilde{\mu}(\widetilde{X}_n)=1$, such that (\ref{e4-11}) and (\ref{e4-12}) hold for $\{(x,u,v_n):(x,u)\in\widetilde{X}_n\}$. Set
    $$\widetilde{X}_0=\bigcap\widetilde{X}_n$$
    Then $\widetilde{\mu}(\widetilde{X}_0)=1$. We shall show that (\ref{e-413}) holds for every $(x,u)\in \widetilde{X}_0$.

Fix $(x,u)\in\widetilde{X}_0$. For every $n,N\in\mathbb{N}$ large enough such that
$\frac{2\beta}{n}<\min\{a,1-a\}$, considering (\ref{e4-11}), one has
\begin{align}\label{e-4-14}
	&\frac{1}{N}\int _0^N1_{B\times [\frac{2\beta}{n},a]\times (\frac{1}{n},\frac{2}{n})}\big(\Phi_{s\beta,s}(x,u,v_n)\big)ds\nonumber\\
	=&\frac{1}{N}\int_0^N 1_B(T^{(\lfloor s+v_n\rfloor,\lfloor s\beta+u\rfloor)}x)
1_{[\frac{2\beta}{n},a]}(\{s\b+u\}) 1_{(\frac{1}{n},\frac{2}{n})}(\{s+v_n\})ds\nonumber\\
	=&\frac{1}{N}\sum_{k=0}^{N-1}\int_{k+(1/n)-v_n}^{k+(2/n)-v_n}1_B(T^{(\lfloor s+v_n\rfloor,\lfloor s\beta+u\rfloor)}x)1_{[\frac{2\beta}{n},a]} (\{s\beta+u\})ds.
\end{align}
Considering (\ref{e4-12}), one has
\begin{align}\label{e-4-14-2}
	&\frac{1}{N}\int _0^N1_{B\times [0,a+\frac{2\beta}{n}]\times (\frac{1}{n},\frac{2}{n})}\big(\Phi_{s\beta,s}(x,u,v_n)\big)ds\nonumber\\
	=&\frac{1}{N}\sum_{k=0}^{N-1}\int_{k+(1/n)-v_n}^{k+(2/n)-v_n}1_B(T^{(\lfloor s+v_n\rfloor,\lfloor s\beta+u\rfloor)}x)1_{[0,a+\frac{2\beta}{n}]} (\{s\beta+u\})ds.
\end{align}
Since $v_n\in(0,1/n^2)$, for every $k\in\{0,1,2,\cdots,N-1\}$ and $s\in [k+\frac{1}{n}-v_n,k+\frac{2}{n}-v_n]$,
we have $(s-k)\beta\in [(\frac{1}{n}-v_n)\beta,(\frac{2}{n}-v_n)\beta]\subset [0,\frac{2\beta}{n}]$. We note the following observation
\begin{itemize}
	\item [(i)] $\lfloor s+v_n\rfloor=k$,
	\item [(ii)] $1_{[\frac{2\beta}{n},a]}(\{s\beta+u\})=1$ implies $1_{[0,a]}(\{k\beta+u\})=1$, and in this case, $\lfloor s\beta+u\rfloor=\lfloor k\beta+u\rfloor$,
	\item [(iii)] $1_{[0,a]}(\{k\beta+u\})=1$ implies $1_{[0,a+\frac{2\beta}{n}]}(\{s\beta+u\})=1$, and in this case, $\lfloor s\beta+u\rfloor=\lfloor k\beta+u\rfloor$.
\end{itemize}
Hence,
\begin{align}\label{e-4-16}
	(\ref{e-4-14})&\overset{\text{(i),(ii)}}{\leq}\frac{1}{N}\sum_{k=0}^{N-1}\int_{k+(1/n)-v_n}^{k+(2/n)-v_n}1_B(T^{(k,\lfloor k\beta+u\rfloor)}x)1_{[0,a]} (\{k\beta+u\})ds \nonumber\\
	&=\frac{1}{N}\sum_{k=0}^{N-1}\frac{1}{n}\cdot 1_B(T^{(k,\lfloor k\beta+u\rfloor)}x)1_{[0,a]} (\{k\beta+u\})\nonumber\\
	&=\frac{1}{n}\cdot \frac{1}{N}\sum_{k=0}^{N-1}1_{B\times [0,a]}(\widetilde{T}^{k}(x,u)).
\end{align}

Similarly with (\ref{e-4-14}) and (\ref{e-4-16}), considering (\ref{e-4-14-2}) we also have
\begin{align}\label{e-4-15}
   (\ref{e-4-14-2})\overset{\text{(i),(iii)}}{\geq}& \frac{1}{N}\sum_{k=0}^{N-1}\int_{k+(1/n)-v_n}^{k+(2/n)-v_n}1_B(T^{(k,\lfloor k\beta+u\rfloor)}x)1_{[0,a]} (\{k\beta+u\})ds\nonumber\\
   =&\frac{1}{n}\cdot \frac{1}{N}\sum_{k=0}^{N-1}1_{B\times [0,a]}(\widetilde{T}^{k}(x,u)).
\end{align}

Combining with (\ref{e4-11}), (\ref{e4-12}) and (\ref{e-4-16}), (\ref{e-4-15}), it is easy to obtain that (\ref{e-413}) holds for every $(x,u)\in\widetilde{X}_0$. This ends our proof.
\end{proof}

In \cite[Theorem 5.1]{RRS2023} Robinson, Rosenblatt and \c{S}ahin proved that for ergodic $\mathbb{Z}^2$-MPS $(X,\mathcal{B}_X,\mu,T)$, there are
at most countably many directions for which $(X,\mathcal{B}_X,\mu,T)$ is not directional ergodic in terms of spectrum.
Therefore, by Theorem \ref{thm:ve=>ves} there are
at most countably many directions such that $(X,\mathcal{B}_X,\mu,T)$ is not directional ergodic.

\begin{thm}\cite[Theorem 5.1]{RRS2023}\label{thm:at most no ergodic}
    Let  $(X,\mathcal{B}_X,\mu,T)$ be an ergodic $\mathbb{Z}^2$-MPS.
Then there are at most countably many directions $\vec{v}$ such that $(X,\mathcal{B}_X,\mu,T)$ is not $\vec{v}$-ergodic.
\end{thm}


\subsection{Directional weak mixing}
In this subsection, we prove the equivalence between directional weak mixing and directional weak mixing in terms of spectrum.
\begin{thm}\label{thm:weak mixing via spectrum}
      Let $(X,\mathcal{B}_X,\mu,T)$ be a $\mathbb{Z}^2$-MPS and $\vec{v}=(1,\beta)$ with $\beta \not\in \mathbb{Q}$.
      Then the following two statements are equivalent:
      \begin{enumerate}
          \item [(1)]$(X,\mathcal{B}_X,\mu,T)$ is  $\vec{v}$-weakly mixing;
          \item [(2)]$(X,\mathcal{B}_X,\mu,T)$ is  $\vec{v}$-weakly mixing in terms of spectrum.
      \end{enumerate}
\end{thm}
\begin{proof}
By Theorem~\ref{thm:weak mixing+directional ergodic=directional weak},
$(X,\mathcal{B}_X,\mu,T)$ is $\vec v$-weakly mixing if and only if it is
weakly mixing and $\vec v$-ergodic.
In \cite[Theorem 6.1]{RRS2023}, the authors proved that
$(X,\mathcal{B}_X,\mu,T)$ is $\vec v$-weakly mixing in terms of spectrum if and only if it is
weakly mixing and $\vec v$-ergodic in terms of spectrum.

By Theorem~\ref{thm:ve=>ves},
$\vec v$-ergodicity is equivalent to $\vec v$-ergodicity in terms of
spectrum, thus $\vec v$-weak mixing is equivalent to $\vec v$-weak mixing in terms of
spectrum.

\end{proof}

In \cite[Theorem 5.1]{RRS2023} Robinson, Rosenblatt and \c{S}ahin proved that for a weakly mixing $\mathbb{Z}^2$-MPS $(X,\mathcal{B}_X,\mu,T)$, there are
at most countably many directions for which $(X,\mathcal{B}_X,\mu,T)$ is not directional weak mixing in terms of spectrum.
Therefore, by Theorem \ref{thm:weak mixing via spectrum} there are
at most countably many directions such that $(X,\mathcal{B}_X,\mu,T)$ is not directional weak mixing.

\begin{thm}\cite[Theorem 5.1]{RRS2023}\label{thm:at most no weak mixing}
    Let  $(X,\mathcal{B}_X,\mu,T)$ be a weakly mixing $\mathbb{Z}^2$-MPS.
    Then at most countably many $\vec{v}$ satisfy that $(X,\mathcal{B}_X,\mu,T)$ is not $\vec{v}$-weakly mixing.
\end{thm}

\begin{cor}
     Let $(X,\mathcal{B}_X,\mu,T)$ be a weakly mixing $\mathbb{Z}^2$-MPS
     and $\vec{v}=(1,\beta)$ with $\beta \not\in \mathbb{Q}$.
     Then the following four statements are equivalent:
      \begin{enumerate}
          \item $(X,\mathcal{B}_X,\mu,T)$ is  $\vec{v}$-weakly mixing;
          \item $(X,\mathcal{B}_X,\mu,T)$ is  $\vec{v}$-weakly mixing in terms of spectrum;
          \item $(X,\mathcal{B}_X,\mu,T)$ is  $\vec{v}$-ergodic;
          \item $(X,\mathcal{B}_X,\mu,T)$ is  $\vec{v}$-ergodic in terms of spectrum.
      \end{enumerate}
\end{cor}
\section{Some examples}\label{Sec:examples}
The following example shows that the limit function $f^*_b$ may depend on the choice of $b\in(0,\infty)$. In particular, this explains why in the definition of directional ergodicity, we cannot replace ``for every $b>0$" by ``there exists $b>0$".

\begin{exam}\label{ex2}
Fix $\beta \not \in \mathbb{Q}$ and let $\vec v=(1,\beta)$. We give an
example on the circle $\mathbb T=\mathbb R/\mathbb Z$ showing that the
pointwise limit $f_b^*$ may depend on the strip width $b$, even when the
underlying $\mathbb Z^2$-action is ergodic.

Let $m_{\mathbb T}$ denote normalized Lebesgue measure and consider the
irrational circle rotation
\[
 R_\beta\colon\mathbb T\longrightarrow\mathbb T,
 \qquad R_\beta(x)=x+\beta\pmod 1.
\]
Define the measure-preserving $\mathbb Z^2$-action on
$(\mathbb T,m_{\mathbb T})$ by
\[
 T^{(m,n)}x=R_\beta^m x=x+m\beta\pmod 1,
 \qquad (m,n)\in\mathbb Z^2.
\]
This action is ergodic because the irrational rotation $R_\beta$ is
ergodic. Identify $\mathbb T$ with $[0,1)$ when specifying intervals, and
put
\[
 A=\left(\frac15,\frac45\right),\qquad
 B=\left[0,\frac15\right]\cup\left[\frac45,1\right),
 \qquad f= 1_A\in L^1(m_{\mathbb T}).
\]

We first take $b=1$. For each integer $m$, there is exactly one integer
$n$ with $|n-m\beta|\le 1/2$. Since $T^{(m,n)}=R_\beta^m$, we have
\[
 \mathbb E(f,\Lambda_M^{\vec v}(1))(x)
 =\frac1M\sum_{m=0}^{M-1}1_A(x+m\beta).
\]
By Weyl equidistribution, for every $x\in\mathbb T$,
\begin{equation}\label{ex2:wide-limit}
 f_1^*(x)
 =\lim_{M\to\infty}\mathbb E(f,\Lambda_M^{\vec v}(1))(x)
 =m_{\mathbb T}(A)=\frac35.
\end{equation}

Next take $b=2/5$. Since the interval
$[m\beta-1/5,m\beta+1/5]$ has length $2/5<1$, it contains an integer
if and only if $\{m\beta\}\in B$, and that integer is unique. Thus
\begin{equation}\label{ex2:narrow-average}
 \mathbb E\left(f,\Lambda_M^{\vec v}\left(\frac25\right)\right)(x)
 =\frac{\displaystyle\sum_{\substack{0\le m<M\\\{m\beta\}\in B}}
          \mathbf1_A(x+m\beta)}
        {\displaystyle\#\{0\le m<M:\{m\beta\}\in B\}}.
\end{equation}
By Weyl equidistribution, the denominator divided by $M$ converges to
$m_{\mathbb T}(B)=2/5$, while the numerator divided by $M$ converges to
$m_{\mathbb T}(B\cap(A-x))$. A direct computation gives
\[
 m_{\mathbb T}(B\cap(A-x))=
 \begin{cases}
  x,&0\le x\le 2/5,\\
  2/5,&2/5\le x\le 3/5,\\
  1-x,&3/5\le x<1.
 \end{cases}
\]
Consequently, for every $x\in\mathbb T$,
\begin{equation}\label{ex2:narrow-limit}
 f_{2/5}^*(x)
 =\lim_{M\to\infty}
    \mathbb E\left(f,\Lambda_M^{\vec v}\left(\frac25\right)\right)(x)
 =\begin{cases}
    \frac52 x,&0\le x\le 2/5,\\
    1,&2/5\le x\le 3/5,\\
    \frac52(1-x),&3/5\le x<1.
  \end{cases}
\end{equation}

Hence $f_1^*\ne f_{2/5}^*$ in $L^1(m_{\mathbb T})$.
Therefore ``for every $b>0$'' in the definition of directional
ergodicity cannot be weakened to ``there exists $b>0$''.
\end{exam}

The following example exhibits a $\mathbb{Z}^2$-action that is weakly mixing, but fails to be
ergodic along a certain irrational direction.

  \begin{exam}\label{not-derg-Z^2-wm}
  Let $(X,\mathcal{B}_X, \mu, (\Phi_s)_{s\in \mathbb{R}})$ be an $\mathbb{R}$-system satisfying the following conditions:
  \begin{itemize}
  	\item [(i)] it is mixing, that is, for each pair of measurable sets $D$ and $E$, one has $$\lim_{t\to\infty}\mu(\Phi_t(D)\cap E)=\mu(D)\mu(E),$$
  	\item [(ii)] there exists $A_0\in \mathcal{B}_X$ with $\mu(A_0)>0$ such that $\mu \big(\bigcup_{s\in [-1,1]}\Phi_{s}A_0\big)<1$
  \end{itemize}
  Specific examples of such systems and sets can be found in \cite[Page 10]{AvilaForniUlcigrai2011}. Fix $\beta \in \mathbb{R}\setminus \mathbb{Q}$ and let $(X,\mathcal{B}_X, \mu, T)$ be a $\mathbb{Z}^2$-action with
  $$T^{(m,n)}(x)\colon=\Phi_{m\beta-n}(x)\ \ \text{for every}\ \ x\in X \ \text{and}\ \ (m,n)\in\mathbb{Z}^2.$$
Let $\vec{v}=(1,\beta)$. We show that $(X,\mathcal{B}_X, \mu, T)$ is $\mathbb{Z}^2$-weakly mixing but not $\vec{v}$-ergodic, hence not $\vec{v}$-weakly mixing.

We first show $(X,\mathcal{B}_X, \mu, T)$ is not $\vec{v}$-ergodic.
In fact, by Theorem \ref{cor:1}, we only need to show that the corresponding skew product system $(\widetilde{X},\widetilde{\mathcal{B}},\widetilde{\mu},\widetilde{T})$ is not ergodic. Let $\widetilde{A}=A_0\times \mathbb{T}$. Then $\widetilde{\mu}(\widetilde{A})>0$. For every $m\in\mathbb{Z}$, we have
  $$\widetilde{T}^m(x,u)=(T^{(m,\lfloor m\beta+u\rfloor)}x,\{m\beta +u\})=(\Phi_{m\beta-\lfloor m\beta+u\rfloor}x,\{m\beta +u\}).$$
  Note that
  $$-1< m\beta-(m\beta+u)\leq m\beta-\lfloor m\beta+u\rfloor\leq m\beta-(m\beta+u-1)\leq 1.$$
  Thus $\widetilde{T}^m\widetilde{A} \subset \bigcup_{s\in [-1,1]}\Phi_{s}(A_0)\times \mathbb{T} $. Then $$\widetilde{\mu}\big(\bigcup_{m\in \mathbb{Z}}\widetilde{T}^m\widetilde{A}\big)\leq\ \mu\big(\bigcup_{s\in [-1,1]}\Phi_{s}A_0\big)<1.$$
  Hence, $(\widetilde{X},\widetilde{\mathcal{B}},\widetilde{\mu},\widetilde{T})$ is not ergodic.

  Now we shall show that $(X,\mathcal{B}_X, \mu, T)$ is a $\mathbb{Z}^2$-weakly mixing system. For every $A,B\in \mathcal{B}_X$ and $\varepsilon>0$, since $(X,\Phi_s)_{s\in \mathbb{R}}$ is a mixing $\mathbb{R}$-system, there exists $c>0$ such that for every $|s|\ge c$, $|\mu(\Phi_sA\cap B)-\mu(A)\mu(B)|<\varepsilon$. Note that, for $(m,n)\in  \mathbb{Z}^2$ one has $T^{-(m,n)}A=\Phi_{n-m\beta}A$. Then
\begin{align*}
&\limsup_{M\to \infty}\frac{1}{M^2}\sum_{(m,n)\in[0,M-1]^2\cap \mathbb{Z}^2} \big|\mu (T^{-(m,n)}A\cap B)-\mu(A)\mu(B)\big| \\
\le &\limsup_{M\to \infty}\frac{1}{M^2}\sum_{(m,n)\in[0,M-1]^2\cap \mathbb{Z}^2, |n-m\beta|< c}  \big|\mu (\Phi_{n-m\beta}A\cap B)-\mu(A)\mu(B)\big|\\
&\ \ \ +\limsup_{M\to \infty}\frac{1}{M^2}\sum_{(m,n)\in[0,M-1]^2\cap \mathbb{Z}^2, |n-m\beta|\ge c} \big|\mu (\Phi_{n-m\beta}A\cap B)-\mu(A)\mu(B)\big|\\
\leq & \bigg(\limsup_{M\to \infty}\frac{1}{M^2}\sum_{(m,n)\in[0,M-1]^2\cap \mathbb{Z}^2, |n-m\beta|< c} 2\bigg)+\varepsilon \\
\leq &\limsup_{M\to \infty}\frac{2(2c+1)M}{M^2}+\varepsilon=\varepsilon.
\end{align*}
Since $\varepsilon$ is arbitrary, $\lim_{M\to \infty}\frac{1}{M^2}\sum_{(m,n)\in[0,M-1]^2\cap \mathbb{Z}^2} \big|\mu (T^{-(m,n)}A\cap B)-\mu(A)\mu(B)\big|=0$,
which means $(X,\mathcal{B}_X,\mu,T)$ is a $\mathbb{Z}^2$-weakly mixing system.
\end{exam}

\section{Proof of Theorem~\ref{main-thm2}}\label{sec:main-thm2}
In this section, we present the proof of Theorem~\ref{main-thm2}, mainly with the aid of the spectral theory. In fact, we shall establish a stronger statement
(Theorem~\ref{thm:zero-width-limit}).

Let $(X,\mathcal{B}_X,\mu,T)$ be a $\mathbb{Z}^2$-MPS, for every $(m,n)\in\mathbb{Z}^2$, the Koopman operator is defined by $U^{(m,n)}f=f\circ T^{(m,n)}$.
Then there exists a unique joint projection-valued spectral measure $\mathsf P$ of this Koopman representation (\cite[Theorem 5.21]{Sch12}). That is,
$$\mathsf P\colon \mathcal{B}_{\mathbb{T}^2}\to \text{Proj}(L^2(\mu)),$$
satisfying
\begin{equation}\label{eq:joint-spectral-representation}
 U^{(m,n)}
 =\int_{\mathbb{T}^2}e^{2\pi i(m\xi+n\eta)}\mathrm{d}\mathsf P(\xi,\eta), \ \text{for every}\ \ (m,n)\in\mathbb{Z}^2,
\end{equation}
where $\text{Proj}(L^2(\mu))$ is the collection of orthogonal projection operators on $L^2(\mu)$. For every $f\in L^2(\mu)$, we can define its associated scalar spectral measure $\sigma_f$ by
\begin{align}\label{asso-spec-measure-projection}
	\sigma_f(B):=\ip{\mathsf P(B)f}{f},
 \qquad B\in \mathcal{B}_{\mathbb{T}^2}.
\end{align}
Moreover, for any bounded measurable functions $\phi,\psi\colon\mathbb{T}^2\to\mathbb{C}$ one has
\begin{align}\label{int-int2}
\left(\int_{\mathbb T^2}\phi\,\mathrm d\mathsf P\right)
 \left(\int_{\mathbb T^2}\psi\,\mathrm d\mathsf P\right)
 =\int_{\mathbb T^2}\phi\psi\,\mathrm d\mathsf P.
\end{align}
For further details on this part, the reader is referred to \cite{Sch12}.


We shall use the following elementary weighted form of Weyl's theorem.
If $G:\mathbb{T}\to\mathbb C$ is Riemann integrable, the $k$-th Fourier coefficient of $G$ is
\[
 \widehat G(k):=\int_0^1G(x)e^{-2\pi ikx} \mathrm{d} x.
\]
Then for every $\alpha\in\R$ and $\beta\in\mathbb{R}\setminus \mathbb{Q}$,
\begin{equation}\label{eq:weighted-weyl}
 \lim_{M\to\infty}\frac1M\sum_{m=0}^{M-1}
 e^{2\pi im\alpha}G(m\beta)
 =
 \begin{cases}
  \widehat G(k),
   &\alpha+k\beta\in \mathbb{Z}\text{ for some }k\in \mathbb{Z},\\
  0,&\text{otherwise}.
 \end{cases}
\end{equation}
The integer $k$, if it exists, is unique.  
Finally we define
\begin{equation}\label{eq:resonant-parametrization}
 \gamma_{\vec v}(r):=(-\beta r,r)+\mathbb{Z}^2,
 \qquad r\in\R.\nonumber
\end{equation}
$\gamma_{\vec v}$ is injective since $\beta \not \in \mathbb{Q}$. Then we set
\begin{equation}\label{eq:resonant-set}
 \Gamma_{\vec v}:=\gamma_{\vec v}(\R)\in\mathcal{B}_{\mathbb{T}^2}.\nonumber
\end{equation}
We shall show that $\mathsf P(\Gamma_{\vec{v}})$ is actually a conditional expectation with some sub $\sigma$-algebra.

We recall that, if $\mathcal A\subseteq\mathcal B_X$ is a sub-$\sigma$-algebra, the conditional expectation $\mathbb E_\mu(\cdot\mid\mathcal A)$ is characterized by $\mathcal A$-measurability and
\[
 \int_A\mathbb E_\mu(f\mid\mathcal A)\,\mathrm d\mu=\int_Af\,\mathrm d\mu,
 \qquad A\in\mathcal A
\]
for every $f\in L^1(\mu)$.

\begin{thm}\label{thm:zero-width-limit} Let $(X,\mathcal{B}_X,\mu,T)$ be a $\mathbb{Z}^2$-MPS, and $\vec{v}=(1,\beta)\in\mathbb{R}^2$ be a direction vector with $\beta \not\in \mathbb{Q}$.
There exists a $T$-invariant sub-$\sigma$-algebra $\mathcal{F}_{\vec v}\subseteq \mathcal{B}_X$, unique modulo $\mu$-null sets, such that for every $f\in L^1(\mu)$,
\begin{equation}\label{eq:zero-width-L1}
 \lim_{b\to0^+}
 \norm{f_b^*-\mathbb E_\mu(f\mid \mathcal{F}_{\vec v})}_1=0.
\end{equation}
On $L^2(\mu)$, this conditional expectation is the spectral projection
onto $\Gamma_{\vec v}$:
\begin{equation}\label{eq:conditional-expectation-spectral-projection}
 \mathbb E_\mu(\,\cdot\mid \mathcal{F}_{\vec v})
 =\mathsf P(\Gamma_{\vec v}).
\end{equation}
Consequently, $\mathcal{F}_{\vec v}$ is trivial if and only if $(X,\mathcal{B}_X,\mu,T)$ is $\vec v$-ergodic.
\end{thm}

\begin{proof} Let $\mathsf P_{\vec{v}}$ denote the orthogonal spectral projection $\mathsf P(\Gamma_{\vec{v}})$.
For $f\in L^1(\mu)$, let $f^*_b$ be the limit function in Theorem \ref{thm-1}. We divide the proof into four steps.

\medskip
\noindent\textbf{Step 1:} In this step, we shall show that for every $f\in L^2(\mu)$, one has
\begin{align}\label{l2-convergence}
	\Vert f^*_b-\mathsf P_{\vec{v}}f\Vert_2\to 0 \ \ \text{as}\ b\to 0.
\end{align}

For $0<b<1$ and $M\in\mathbb{N}$, we consider the operator on $L^2(\mu)$
\begin{equation}\label{eq:finite-strip-average}
 \mathcal A_{M,b}f
 :=\mathbb E(f,\Lambda_M^{\vec v}(b))
 =\frac1{\#\Lambda_M^{\vec v}(b)}
   \sum_{(m,n)\in\Lambda_M^{\vec v}(b)}U_{T^{(m,n)}}f.
\end{equation}
Formula \eqref{eq:joint-spectral-representation} shows that
\begin{align}\label{A_n-a_n}
	\mathcal{A}_{M,b}=\int_{\mathbb{T}^2} a_{M,b}(z)\mathrm{d} \mathsf P(z),
\end{align}
where \begin{equation}\label{eq:scalar-multiplier}
 a_{M,b}(z)
 =\frac{1}{\#\Lambda_M^{\vec v}(b)}\cdot
   \sum_{(m,n)\in\Lambda_M^{\vec v}(b)}
   e^{2\pi i(m\xi+n\eta)}\ \ \text{for every}\ \ z=(\xi,\eta)\in\mathbb{T}^2.\nonumber
\end{equation}
For every $b>0$ and $\eta\in [0,1)$, we define \begin{equation}\label{eq:auxiliary-periodic-function}
 G_{b,\eta}(x)
 :=\sum_{n\in \mathbb{Z}}
   1_{[-b/2,b/2]}(n-x)e^{2\pi i\eta(n-x)},
 \qquad x\in \mathbb{T}.\nonumber
\end{equation}
It is well-defined, being independent of the choice of representative of $x$. For each $x$ the sum is finite. Moreover, $G_{b,\eta}$ is bounded and Riemann integrable. Then for every $M\in\mathbb{N}$, $b>0$ and $z=(\xi,\eta)\in\mathbb{T}^2$
\begin{align*}
	a_{M,b}(z)=&\frac{1}{\#\Lambda_M^{\vec v}(b)}\sum_{m=0}^{M-1}\sum_{n\in\mathbb{Z}}1_{[-b/2,b/2]}(n-m\beta)\cdot e^{2\pi i(m\xi+n\eta)}\\
	=&\frac{1}{\#\Lambda_M^{\vec v}(b)}\sum_{m=0}^{M-1}\sum_{n\in\mathbb{Z}}1_{[-b/2,b/2]}(n-m\beta)\cdot e^{2\pi i \eta(n-m\beta)}\cdot e^{2\pi im(\xi+\eta\beta)}\\
	=&\frac{1}{\#\Lambda_M^{\vec v}(b)}\sum_{m=0}^{M-1}G_{b,\eta}(m\beta)\cdot e^{2\pi im(\xi+\eta\beta)}.
\end{align*}
Then by (\ref{eq:weighted-weyl}), the limit of $a_{M,b}(z)$ as $M$ goes to infinity can be non-zero only if   $(\xi+\eta\beta)+k\beta\in\mathbb{Z}$ for some $k\in\mathbb{Z}$. Equivalently,
$$z=(\xi,\eta)+\mathbb{Z}^2=(-\beta(\eta+k),\eta+k)+\mathbb{Z}^2=\gamma_{\vec v}(\eta+k)\in\Gamma_{\vec v}.$$

The $k$-th Fourier coefficient of $G_{b,\eta}$ is 
\begin{align}
 \widehat G_{b,\eta}(k)
 &=\int_0^1\sum_{j\in \mathbb{Z}}
   1_{[-b/2,b/2]}(j-x)e^{2\pi i\eta(j-x)}
   e^{-2\pi ikx}\mathrm{d} x \notag\\
   &=\sum_{j\in \mathbb{Z}}\int_0^1
   1_{[-b/2,b/2]}(j-x)e^{2\pi i\eta(j-x)}
   e^{-2\pi ikx}\mathrm{d} x \notag\\
   &=\sum_{j\in \mathbb{Z}}\int_{j-1}^j
   1_{[-b/2,b/2]}(t)e^{2\pi i\eta t}
   e^{2\pi ik(t-j)}\mathrm{d} t \notag\\
   &=\sum_{j\in \mathbb{Z}}\int_{j-1}^j
   1_{[-b/2,b/2]}(t)e^{2\pi i\eta t}
   e^{2\pi ikt}\mathrm{d} t \notag\\
 &=\int_{\mathbb{R}}
   1_{[-b/2,b/2]}(t)e^{2\pi i\eta t}
   e^{2\pi ikt}\mathrm{d} t \notag
   =\int_{-b/2}^{b/2}e^{2\pi i(\eta+k)t}\mathrm{d} t.
\end{align}
Thus $a_{M,b}(z)\to\varphi_b(z)$ as $M\to\infty$, where
\begin{equation}\label{eq:phi-b-definition}
\varphi_b(z)
 :=
 \begin{cases}
  \displaystyle
  \frac1b\int_{-b/2}^{b/2}e^{2\pi irt}\mathrm{d} t
  =\frac{\sin(\pi br)}{\pi br},
  &z=\gamma_{\vec v}(r),\quad r\neq0,\\[3mm]
  1, &z=(0,0),\\
  0,&z\notin\Gamma_{\vec v}.
 \end{cases}
\end{equation}
Since $\abs{a_{M,b}}\leq1$ and $\abs{\varphi_b}\leq1$,
by  \eqref{asso-spec-measure-projection},  \eqref{int-int2}, \eqref{A_n-a_n} and dominated convergence theorem,
for every
$f\in L^2(\mu)$,
\begin{align}
 \norm{\mathcal{A}_{M,b}f-
  \int_{\mathbb{T}^2}\varphi_b(z)\mathrm{d}\mathsf P(z)f}_2^2
 &=\int_{\mathbb{T}^2}\abs{a_{M,b}-\varphi_b}^2\mathrm{d}\sigma_f
 \longrightarrow0,\ \text{as}\ M\to\infty.
 \label{eq:fixed-width-L2-spectral-limit}\nonumber
\end{align}
On the other hand, by Theorem~\ref{D-von Ne-thm}, $\mathcal{A}_{M,b}f=\mathbb E(f,\Lambda_M^{\vec v}(b))\to f_b^*$ in $L^2(\mu)$. By uniqueness of the $L^2$
limit, one has
\begin{equation}\label{eq:fixed-width-spectral-formula}
 f^*_b
 =\int_{\mathbb{T}^2}\varphi_b(z)\mathrm{d}\mathsf P(z)f,
 \qquad \text{for every}\ \ f\in L^2(\mu).
\end{equation}
Moreover, since $\varphi_b\to\mathbf{1}_{\Gamma_{\vec v}}$ pointwise and
$\abs{\varphi_b}\leq1$, by dominated
convergence theorem and
\eqref{eq:fixed-width-spectral-formula}
\begin{equation}\label{eq:L2-zero-width-limit}
 \norm{f^*_b-\mathsf P_{\vec v}f}_2^2
 =\int_{\mathbb{T}^2}
   \abs{\varphi_b-\mathbf{1}_{\Gamma_{\vec v}}}^2\mathrm{d}\sigma_f
 \longrightarrow0
 \qquad(b\to0^+)
\end{equation}
for every $f\in L^2(\mu)$.

\medskip
\noindent\textbf{Step 2:} We set $\Ran \mathsf P_{\vec v}\colon=\{\mathsf P_{\vec v}f:f\in L^2(\mu)\}$ and $\mathcal{F}_{\vec v}
 :=\bigl\{A\in \mathcal{B}_X:1_A\in\Ran \mathsf P_{\vec v}\bigr\}$. In this step, we shall show that $\mathcal{F}_{\vec{v}}$ is a $T$-invariant sub-$\sigma$-algebra of $\mathcal{B}_X$, $\Ran \mathsf P_{\vec v}=L^2(\mathcal{F}_{\vec v})$ and
 \begin{align}\label{l2-conditional exp}
 	 \mathsf P_{\vec v}=\mathbb E_\mu(\,\cdot\mid\mathcal F_{\vec v})
 \qquad\text{on }L^2(\mu).
 \end{align}

First, $\mathsf P_{\vec v}$ is positive. Indeed, if $h\in L^2(\mu)$ and $h\ge0$, then $h^*_b\ge0$ for every $b>0$, while \eqref{l2-convergence} gives $h^*_b\to\mathsf P_{\vec v}h$ in $L^2$. Then $\mathsf P_{\vec v}h\ge0$.
Similarly, since $1^*_b=1$, we have $\mathsf P_{\vec v}1=1$.

Let $g\in\Ran \mathsf P_{\vec v}$ be real-valued. Positivity gives
\[
 \abs{g}=\abs{\mathsf P_{\vec v}g}\leq \mathsf P_{\vec v}\abs{g}.
\]
Since $\mathsf P_{\vec v}$ is self-adjoint and fixes $1$,
\[
 \int_X \mathsf P_{\vec v}\abs{g}\mathrm{d}\mu
 =\ip{\mathsf P_{\vec v}\abs{g}}{1}
 =\ip{\abs{g}}{\mathsf P_{\vec v}1}
 =\int_X\abs{g}\mathrm{d}\mu.
\]
Thus the nonnegative function
$\mathsf P_{\vec v}\abs{g}-\abs{g}$ has integral zero, and consequently
$\mathsf P_{\vec v}\abs{g}=\abs{g}$.  Therefore $\Ran \mathsf P_{\vec v}$ is closed
under absolute values and hence under
$\min$ and $\max$. Then by \[
 1_{A^c}=1-1_A,
 \qquad
 1_{A\cap B}=\min\{1_A,1_B\}, \ \text{for every}\ A,B\in\mathcal{B}_{X},
\]
$\mathcal{F}_{\vec{v}}$ is closed under complements and finite intersections.

If $(A_j)_{j\ge1}\subseteq\mathcal F_{\vec v}$ are pairwise disjoint and $B_N=\bigcup_{j=1}^NA_j$, then $1_{B_N}\in\operatorname{Ran}\mathsf P_{\vec v}$ and $1_{B_N}\to1_{\cup_jA_j}$ in $L^2$. Since the range of an orthogonal projection is closed, $\cup_jA_j\in\mathcal F_{\vec v}$. Hence $\mathcal F_{\vec v}$ is a $\sigma$-algebra.
Simple functions and $L^2$ approximation show that $L^2(\mathcal{F}_{\vec v})\subseteq\Ran \mathsf P_{\vec v}$.

Moreover, by \eqref{eq:joint-spectral-representation} and
\eqref{int-int2}, for every $(m,n)\in\mathbb{Z}^2$ we have
\[
 \mathsf P_{\vec v}U_{T^{(m,n)}}
 =U_{T^{(m,n)}}\mathsf P_{\vec v}.
\]
Let $A\in\mathcal F_{\vec v}$. Since
$\mathsf P_{\vec v}1_A=1_A$ and
$U_{T^{(m,n)}}1_A=1_{(T^{(m,n)})^{-1}A}$, it follows that
\[
 \mathsf P_{\vec v}1_{(T^{(m,n)})^{-1}A}
 =\mathsf P_{\vec v}U_{T^{(m,n)}}1_A
 =U_{T^{(m,n)}}\mathsf P_{\vec v}1_A
 =1_{(T^{(m,n)})^{-1}A}.
\]
Thus $(T^{(m,n)})^{-1}A\in\mathcal F_{\vec v}$ for all
$(m,n)\in\mathbb{Z}^2$. Thus $\mathcal F_{\vec v}$ is $T$-invariant.

Conversely,
let $g\in\operatorname{Ran}\mathsf P_{\vec v}$ be real-valued. For $q\in\mathbb Q$ and $k\ge1$,
since $\Ran \mathsf P_{\vec v}$ is closed
under $\min$ and $\max$,
\[
 h_{k,q}:=\min\{1,\max\{0,k(g-q)\}\}\in\operatorname{Ran}\mathsf P_{\vec v}.
\]
Since $h_{k,q}\to1_{\{g>q\}}$ in $L^2(\mu)$, we have $\{g>q\}\in\mathcal F_{\vec v}$. Thus $g$ is $\mathcal F_{\vec v}$-measurable. The complex case follows by considering real and imaginary parts,
thus $L^2(\mathcal{F}_{\vec v})= \Ran \mathsf P_{\vec v}$.

Finally, if $D\in\mathcal F_{\vec v}$ and $f\in L^2(\mu)$, then $\mathsf P_{\vec v}1_D=1_D$, and self-adjointness gives
\[
 \int_D\mathsf P_{\vec v}f\,\mathrm d\mu
 =\ip{\mathsf P_{\vec v}f}{1_D}
 =\ip{f}{\mathsf P_{\vec v}1_D}
 =\int_Df\,\mathrm d\mu.
\]
Hence
$\mathsf P_{\vec v}=\mathbb E_\mu(\,\cdot\mid\mathcal F_{\vec v})$.

\medskip
\noindent\textbf{Step 3:} In this step, we will show that \eqref{eq:zero-width-L1} holds for every $f\in L^1(\mu)$.
Let $f\in L^1(\mu)$ and let $\varepsilon>0$.  Choose
$g\in L^2(\mu)$ such that $\norm{f-g}_1<\varepsilon$.  By Proposition~\ref{continuous-on-f}
\begin{equation}
 \norm{f^*_b-g^*_b}_1\leq\norm{f-g}_1.\nonumber
\end{equation}
Conditional expectation is also an $L^1$ contraction. Hence, using
(\ref{l2-conditional exp}),
\begin{align}
 \norm{f^*_b-\mathbb E_\mu(f\mid \mathcal{F}_{\vec v})}_1
 &\leq
 \norm{f^*_b-g^*_b}_1
 +\norm{g^*_b-\mathsf P_{\vec v}g}_1 \notag
 +\norm{\mathsf P_{\vec v}g-
        \mathbb E_\mu(f\mid \mathcal{F}_{\vec v})}_1 \notag\\
        &=
 \norm{f^*_b-g^*_b}_1
 +\norm{g^*_b-\mathsf P_{\vec v}g}_1 \notag
 +\norm{\mathbb E_\mu(g\mid \mathcal{F}_{\vec v})-
        \mathbb E_\mu(f\mid \mathcal{F}_{\vec v})}_1 \notag\\
 &\leq 2\norm{f-g}_1+\norm{g^*_b-\mathsf P_{\vec v}g}_1 \notag
 \leq 2\varepsilon+\norm{g^*_b-\mathsf P_{\vec v}g}_2.\nonumber
\end{align}
By \eqref{l2-convergence},
 $\limsup_{b\to0^+}\norm{f^*_b-\mathbb E_\mu(f\mid \mathcal{F}_{\vec v})}_1\le 2\varepsilon$. The conclusion follows from the arbitrariness of $\varepsilon$.

\medskip
\noindent\textbf{Step 4:} In this step, we will show that $(X,\mathcal{B}_X,\mu,T)$ is $\vec v$-ergodic if and only if
$\mathcal{F}_{\vec v}$ is trivial.

We show that for every $b>0$ and $f\in L^2(\mu)$, $f^*_b\in\Ran \mathsf P_{\vec{v}}$. In fact, by (\ref{eq:fixed-width-spectral-formula})
\begin{align}\label{20261024}
	\mathsf P_{\vec{v}}f^*_b=&\mathsf P_{\vec{v}}\int_{\mathbb{T}^2}\varphi_b(z)\mathrm{d}\mathsf P(z)f=\big(
\int_{\mathbb{T}^2}1_{\Gamma_{\vec{v}}}(z) \mathrm{d}\mathsf{P}(z)\big)
\big(\int_{\mathbb{T}^2}\varphi_b(z)\mathrm{d}\mathsf{P}(z)\big)f \nonumber\\
\overset{(\ref{int-int2})}{=}&\int_{\mathbb{T}^2}1_{\Gamma_{\vec{v}}}\varphi_b(z)\mathrm{d}\mathsf{P}(z)f
\overset{(\ref{eq:phi-b-definition})}{=}\int_{\mathbb{T}^2}\varphi_b(z)\mathrm{d}\mathsf{P}(z)f=f^*_b.
\end{align}
This means $f^*_b\in \Ran \mathsf P_{\vec{v}}$ for every $f\in L^2(\mu)$.

If $\mathcal{F}_{\vec{v}}$ is trivial, by (\ref{l2-conditional exp}),
$\Ran \mathsf P_{\vec v}=\{\mathsf P_{\vec v}f:f\in L^2(\mu)\}$ contains only constant functions. Then by \eqref{20261024}, $f^*_b=\int_X fd\mu$ for every $f\in L^2(\mu)$ and $b>0$. Noting that
$L^2(\mu)$ is dense in $L^1(\mu)$,
by Proposition \ref{continuous-on-f}, $f^*_b=\int_X fd\mu$ for every $f\in L^1(\mu)$ and $b>0$.
Thus $(X,\mathcal{B}_X,\mu,T)$ is $\vec v$-ergodic.

If $(X,\mathcal{B}_X,\mu,T)$ is $\vec v$-ergodic,
then $f^*_b=\int_X fd\mu$ for every $f\in L^2(\mu)$ and $b>0$.
By \eqref{eq:L2-zero-width-limit}
 $\mathsf P_{\vec{v}}f=\int_X fd\mu$. By (\ref{l2-conditional exp}), $\mathcal{F}_{\vec{v}}$ is trivial.
\end{proof}

\section{Acknowledgement}
 We are grateful to Professors Wen Huang and Hanfeng Li for their helpful suggestions and insightful comments. 

After completing this paper, we shared our manuscript with Professors E. A. Robinson and A. A. \c{S}ahin, who have conducted extensive research on directional ergodicity and weak mixing \cite{RRS2023}. They kindly shared with us their recent work \cite{RS2026}, which provides valuable insights into directional ergodicity and weak mixing from a spectral-theoretic perspective. In particular, they established the equivalence between dynamical and spectral directional weak mixing for $\mathbb{Z}^2$-actions.
Their proof \cite{RS2026} is based on an approach fundamentally different from the one developed in the present paper.

\section{AI disclosure}
ChatGPT (GPT-5.6 Plus) was used to assist with language polishing, typographical
corrections, and final consistency checks in the whole paper.

The authors developed the initial argument for
Theorem \ref{main-thm2} and used ChatGPT (GPT-5.6 Pro)
for assistance in revising and completing its exposition. The authors
independently checked the mathematical arguments and take full responsibility for the paper.

\end{document}